\documentclass[11pt,a4paper]{article}

\usepackage[title]{appendix}
\usepackage[normalem]{ulem}
\usepackage{graphicx}
\usepackage{amsmath}
\usepackage{amsfonts}
\usepackage{amssymb}
\usepackage{enumerate}
\usepackage{lscape}
\usepackage{longtable}
\usepackage{rotating}
\usepackage{color}
\usepackage{url}
\usepackage{authblk}
\usepackage[font={small,it}]{caption}
\usepackage{subcaption}
\usepackage{rotating}
\usepackage{titlesec}
\usepackage[title]{appendix}
\newcommand{\norm}[1]{\left\lVert#1\right\rVert}
\usepackage{multicol}
\usepackage{mathrsfs}
\usepackage{upgreek}
\usepackage{lmodern}
\usepackage{algorithmicx}
\usepackage{algorithm}
\usepackage[noend]{algpseudocode}

\newtheorem{theorem}{Theorem}[section]

\newtheorem{lemma}{Lemma}[section]

\newtheorem{proposition}{Proposition}[section]

\newtheorem{assumption}{Assumption}[section]

\newenvironment{proof}[1][Proof]{\textbf{#1.} }{\hfill$\Box$}

\makeatletter
\newcommand*{\centerfloat}{%
  \parindent \z@
  \leftskip \z@ \@plus 1fil \@minus \textwidth
  \rightskip\leftskip
  \parfillskip \z@skip}
\makeatother

\usepackage{soul}
\newcommand{\FLE}{\texttt{Full-Low Evaluation}}

\newcommand{\FullEval}{\texttt{Full-Eval}}
\newcommand{\LowEval}{\texttt{Low-Eval}}

\newcommand{\NOMAD}{\texttt{NOMAD}}
\newcommand{\ps}{\texttt{patternsearch}}

\newcommand\cuter{{\sf CUTEst}}

\DeclareMathOperator{\argmin}{argmin}

\def\low{\mathrm{low}}

\def\R{\mathbb{R}}

\def\calF{\mathcal{F}}
\def\T{\mathrm{T}}
\newcommand{\PF}[1]{P_{\mathcal{F}}\left[#1\right]}

\numberwithin{equation}{section}
\definecolor{forestgreen}{rgb}{0.13, 0.55, 0.13}
\definecolor{greenp}{rgb}{0.0, 0.65, 0.31}
\definecolor{greenr}{rgb}{0.4, 0.69, 0.2}
\definecolor{indiagreen}{rgb}{0.07, 0.53, 0.03}
\definecolor{Blue}{rgb}{0, 0, 1}

\usepackage[square,numbers]{natbib}
\usepackage{hyperref}
\hypersetup{
  colorlinks,
  citecolor=indiagreen,
  linkcolor=indiagreen,
  urlcolor=Blue}

\newcommand{\revisedII}[1]{\textcolor{red}{#1}}
\newcommand{\ceil}[1]{\left\lceil #1 \right\rceil}

\begin{document}

\title{Full-low evaluation methods for bound and linearly constrained derivative-free 
optimization}


\author[1]{C. W. Royer}
\author[2]{O. Sohab\thanks{Corresponding author. Email: {\tt ous219@lehigh.edu}}}
\author[2]{L. N. Vicente}

\affil[1]{LAMSADE, CNRS, Universit\'e Paris Dauphine-PSL, Place du Mar\'echal de Lattre de Tassigny, 75016 Paris, France}
\affil[2]{Department of Industrial and Systems Engineering, Lehigh University, 200 West Packer Avenue, Bethlehem, PA 18015-1582, USA}

\maketitle
\renewcommand{\thefootnote}{\fnsymbol{footnote}}
\footnotetext{Contributing authors: {\tt clement.royer@lamsade.dauphine.fr}, {\tt lnv@lehigh.edu}}
\footnotesep=0.4cm

{\small
\begin{abstract}
Derivative-free optimization (DFO) consists in finding the best value of an 
objective function without relying on derivatives. To tackle such problems, 
one may build approximate derivatives, using for instance finite-difference 
estimates. One may also design algorithmic strategies that perform space 
exploration and seek improvement over the current point. The first type of 
strategy often provides good performance on smooth problems but at the expense 
of more function evaluations. The second type is cheaper and 
typically handles non-smoothness or noise in the objective better. Recently, 
full-low evaluation methods have been proposed as a hybrid class of DFO 
algorithms that combine both strategies, respectively denoted as Full-Eval and 
Low-Eval. In the unconstrained case, these methods showed promising numerical 
performance.

In this paper, we extend the full-low evaluation framework to bound and 
linearly constrained derivative-free optimization. We derive convergence 
results for an instance of this framework, that combines finite-difference
quasi-Newton steps with probabilistic direct-search steps. The former are 
projected onto the feasible set, while the latter are defined within tangent cones 
identified by nearby active constraints. We illustrate the practical 
performance of our instance on standard linearly constrained problems, that 
we adapt to introduce noisy evaluations as well as non-smoothness. In all cases, 
our method performs favorably compared to algorithms that rely solely on 
Full-eval or Low-eval iterations.
\end{abstract}
}

\section{Introduction}
\label{sec:intro}

Derivative-Free Optimization (DFO)~\cite{CAudet_WHare_2017,
ARConn_KScheinberg_LNVicente_2009b,ALCustodio_KScheinberg_LNVicente_2017,
JLarson_MMenickelly_SWild_2019,LMRios_NVSahinidis_2013} methods are developed for the minimization of functions whose corresponding derivatives are unavailable for use or expensive to compute. Particularly useful for complex simulation problems, DFO is often employed when the objective function is derived from numerical simulations, making derivatives inaccessible for algorithmic purposes.
The field of derivative-free optimization now spans a wide range of algorithms and has been applied in numerous engineering and applied science fields~\cite{alarie2021two}.
 
In these settings, 
evaluating the objective function represents the main computational 
bottleneck, that must be accounted for while designing DFO algorithms. 
Besides, simulation codes often enforce hard constraints on their parameters, 
typically under the form of bounds or linear relationships, that must be 
satisfied for the simulation to terminate and for function information to 
be obtained. 
As a result, a plurality of DFO schemes have been developed to target constrained problems, with a careful classification of the nature of those constraints. For a comprehensive coverage of constraints in DFO, we refer the reader to existing survey papers~\cite{JLarson_MMenickelly_SWild_2019,SLeDigabel_SMWild_2024}. We review below algorithmic approaches that bear direct relevance to this paper.

Direct-search methods~\cite{TGKolda_RMLewis_VTorczon_2003} are a common choice 
of derivative-free algorithms due to their ease of implementation. These 
iterative methods sample new function evaluations along suitably chosen 
directions at every iteration, in order to find a point at which the function 
value decreases. Direct-search schemes have been endowed with theoretical 
guarantees even in presence of non-smooth objectives, while being 
intrinsically robust to the presence of noise in the function 
evaluations~\cite{CAudet_WHare_2017}. In presence of linear constraints, 
direct-search methods generally use directions that conform to the geometry 
of the feasible set, thereby ensuring feasible descent without relying on derivative information~\cite{kolda2007stationarity}. 
Recent results have proposed probabilistic variants of direct search, in 
which the directions are only guaranteed to be feasible descent with a 
given probability~\cite{gratton2019direct}. A probabilistic direct-search 
iteration can be performed using a significantly smaller number of directions 
(and, thus, of function evaluations) than its deterministic counterpart.

An alternative to direct-search techniques consists in building an approximate 
derivative from function evaluations, which then enables the calculation of 
steps similar to those in the derivative-based setting. Model-based 
derivative-free techniques obey this logic, and rely on  
trust-region globalization arguments from nonlinear optimization to guarantee 
convergence of the methods~\cite{ARConn_KScheinberg_LNVicente_2009b}. As a
result, bounds and linear constraints are typically handled in a similar 
fashion as in the derivative-based 
setting~\cite{ARConn_NIMGould_PhLToint_2000}, even though ad hoc strategies 
have also been considered~\cite{SGratton_PhLToint_ATroltzsch_2011}. Another 
widely common approach consists in using finite differences to approximate 
derivatives, so as to leverage existing algorithms from the derivative-based 
literature~\cite[Chapter 8]{JNocedal_SJWright_2006}. In particular, recent 
advances in applying quasi-Newton updates using finite-difference gradients 
have demonstrated good numerical performance, especially in a smooth 
setting~\cite{ASBerahas_RHByrd_JNocedal_2019,ASBerahas_OSohab_LNVicente_2022,HJMShi_MQXuan_FOztoprak_JNocedal_2023}. This performance is mitigated by the 
inherent cost of finite-difference estimates, that scales at least linearly 
with the dimension, and thus may be expensive to perform in a simulation-based 
environment.

The full-low evaluation (\FLE{}) framework was recently proposed as a 
principled way of combining derivative-free steps with different costs and 
properties~\cite{ASBerahas_OSohab_LNVicente_2022}. In the unconstrained 
setting, it was proposed to instantiate this framework using a BFGS step 
computed through finite differences as well as a probabilistic direct-search 
iteration. This hybrid approach was shown to outperform the individual 
strategies while being competitive with an established solver on smooth and 
piecewise smooth problems.

In this paper, we propose an extension of the \FLE{} framework that handles 
bounds and linear constraints by producing feasible iterates and feasible 
steps. We analyze an instance of this approach that combines projected steps 
built on finite-difference gradient estimates with direct-search 
steps based on probabilistic feasible descent as handled in \cite{gratton2019direct}. The former (\FullEval{} step) is considered expensive in terms of evaluations but provides good performance 
and convergence results in the presence of a smooth objective. The latter 
(\LowEval{} step) is cheaper in terms of evaluations, while being more robust 
to noise or non-smoothness in the objective. Similarly to the unconstrained 
setting, a switching condition determines the nature of the step taken at 
each iteration.

The rest of this paper is organized as follows. Section~\ref{sec:pbtan} states our 
problem of interest, as well as the key geometrical concepts used to design our 
algorithm. Section~\ref{sec:flec} describes our generic \FLE{} framework, as well as 
the two subroutines that define our instance of interest. Section~\ref{sec:cv} provides 
convergence results for both smooth and non-smooth objectives. Section~\ref{sec:numsetup} 
details our implementation and our experimental setup, while the output of our tests is 
analyzed in Section~\ref{sec:numres}. Final remarks are given in Section~\ref{sec:conc}. 
A list of our test problems is provided in Appendix~\ref{sec:pblist}.

\section{Linearly constrained optimization and tangent cones}
\label{sec:pbtan}

The purpose of this section is twofold. First, we describe our problem of interest as 
well as associated optimality measures in Section~\ref{subsec:optimeas}. These concepts 
will serve as a basis for the \FullEval{} part of our algorithm. Secondly, we 
discuss the notion of tangent cones and its connection to feasibility in 
Section~\ref{subsec:tangcone}. Those definitions will be instrumental in designing our
\LowEval{} step based on direct search.

\subsection{Problem and optimality measure}
\label{subsec:optimeas}

In this paper, we are interested in solving linearly constrained problems of the form
\begin{equation}
\label{eq:prob}
\begin{split}
	\displaystyle \min_{x \in \R^n}&~f(x) \\
	\text{s.t.}  &~Ax = b\\
    &\ell \leq A_\mathcal{I} x \leq u,
\end{split}
\end{equation}
where $f: \R^n \to \R$, $A \in \R^{m \times n}$, 
$A_\mathcal{I} \in \R^{m_\mathcal{I} \times n}$, $b \in \R^m$ and 
$(\ell,u) \in  \bar{\R}^{m_\mathcal{I}} \times \bar{\R}^{m_\mathcal{I}}$ where $ \bar{\R} = \R \cup \{-\infty, \infty\}$ with $\ell < u$. 
To encompass bound constrained problems into the general formulation~\eqref{eq:prob}, 
we consider the possibility that the matrix $A$ is empty, in which case $m$ is equal 
to zero. When $m > 0$, the matrix $A$ is assumed to have full row rank, and we let 
$W\in \R^{n \times (n-m)}$ be an orthonormal basis for the null space of $A$. 
When $m=0$, we consider $W$ to be the identity matrix in $\R^n$. Finally, we denote 
the set of feasible points by $\calF$.

Assuming that the function $f$ is continuously differentiable, it is possible to define a 
criticality measure for problem~\eqref{eq:prob} that characterizes first-order 
stationary point. We focus on the quantity $q(\cdot)$ defined by
\begin{equation}
\label{eq:qmeas}
	 \forall x \in \calF, \quad
	 q(x) \; := \; \left\| \PF{x -\nabla f(x)}-x \right\|,
\end{equation}
where $\displaystyle \PF{x} = \argmin_{y \in \calF} \| x - y\|$ is the projection of $x$ onto the feasible region $\calF$ 
and $\|\cdot\|$ denotes the Euclidean norm.
In derivative-free optimization, the metric~\eqref{eq:qmeas} has been employed for 
analyzing the convergence of algorithms designed for the linearly constrained 
setting~\cite{RMLewis_VTorczon_1999, RMLewis_VTorczon_2000}. Although more recent 
approaches have focused on another metric bearing a close connection with the 
direct-search stepsize~\cite{gratton2019direct, TGKolda_RMLewis_VTorczon_2003}, the measure~\eqref{eq:qmeas} is quite common in projected gradient 
techniques~\cite{DPBertsekas_2016}. Since our theory relies on that of projected 
gradient techniques in the smooth setting, we naturally focus on the measure~\eqref{eq:qmeas}.

When the smoothness of the function $f$ is not guaranteed but $f$ is locally Lipschitz 
continuous, necessary optimality condition for problem \eqref{eq:prob} can be 
formulated based on the Clarke-Jahn generalized directional derivative of 
$f$~\cite{jahn1994introduction}. 
For a given point $x \in \mathcal{F}$, a direction $d$ is called feasible at $x$ there exists $\epsilon > 0$ for which 
$x + \epsilon d \in \calF$. Given $x \in \mathcal{F}$
and a feasible 
direction $d$ at $x$, the Clarke-Jahn generalized directional derivative of $f$ 
at $x$ in direction $d$ is defined as 
\begin{equation}
\label{eq:clarke}
f^\circ (x;d) := \limsup_{\begin{array}{c}
    y \to x, y \in \mathcal{F} \\
     t \downarrow 0, y + td \in \mathcal{F}
\end{array} } \frac{f(y+td) - f(y)}{t}.
\end{equation}
Any $x^* \in \mathcal{F}$ such that $f^\circ (x^*;d) \geq 0$ for any feasible 
direction $d$ is called a Clarke-Jahn stationary point.  Note that such a 
condition was recently used in the context of non-smooth optimization with linear
constraints~\cite{beck2020convergence}. In the linearly constrained case, the set of feasible directions at~$x^*$ coincides with the tangent cone~$T(x^*)$. 

\subsection{Approximate Tangent Cones}
\label{subsec:tangcone}

Tangent cones are key concepts to characterize feasibility and optimality in constrained 
optimization~\cite{JNocedal_SJWright_2006}. Although approximate tangent cones have been less 
studied, they have proven quite useful in the context of derivative-free optimization, as they 
characterize directions that are feasible for a given step size, by accounting for constraints 
that are either active or approximately active~\cite{TGKolda_RMLewis_VTorczon_2003}. 
We recall below the key definitions related 
to approximate tangent cones, by following the presentation in Gratton et 
al.~\cite{gratton2019direct}.

For convenience, we will define approximate tangent cones based on a parameterization of the 
feasible set. More precisely, we fix a reference vector $\bar{x} \in \R^{n-m}$ such that 
$A \bar{x} = b$. Then, any feasible point $x \in \calF$ can be written as 
$x=W \tilde{x}+\bar{x}$, where $\tilde{x} \in \R^{n-m}$ is such that
\[
	\ell - A_{\mathcal{I}} \bar{x} 
	\; \leq \; A_{\mathcal{I}} W \tilde{x} 
	\; \le \;  u - A_{\mathcal{I}} \bar{x} .
\]
Using this decomposition, we define the \emph{approximate active inequality constraints} 
at $x = W \tilde{x}+\bar{x}$ according to a step size $\xi>0$ as
\begin{equation}
\label{eq:activecons}
	\Bigg\{
	\begin{aligned}
		I_u(x,\xi) &:=& \left\{ i :~ 
		|u_i -[A_\mathcal{I} \bar{x}]_i - [A_\mathcal{I} W \tilde{x}]_i | 
		\leq \xi \norm{W^\top A_\mathcal{I}^\top e_i} \right\} \\
		I_\ell(x,\xi) &:=& \left\{ i :~
		|\ell_i -[A_\mathcal{I} \bar{x}]_i - [A_\mathcal{I} W \tilde{x}]_i | 
		\leq \xi \norm{W^\top A_\mathcal{I}^\top e_i} \right\},
	\end{aligned}
\end{equation}
where $e_1, \dots, e_{m_\mathcal{I}}$ denote the coordinate vectors in 
$\R^{m_\mathcal{I}}$. Note that one can assume without loss of generality that \(W^\top A_\mathcal{I}^\top e_i \neq 0\), otherwise, given that we assume that \(\mathcal{F}\) is nonempty, the inequality constraints/bounds \(\ell_i \leq [A_\mathcal{I} x]_i \leq u_i\) would be redundant. Those indices in turn define the \emph{approximate normal 
cone} associated with $(x,\xi)$ as
\begin{equation}
\label{eq:normalcone}
	N(x, \xi) \; := \; \text{Cone}\left(
	\{W^\top A^\top_\mathcal{I} e_i \}_{i \in I_u(x, \xi)} \} 
	\cup 
	\{-W^\top A^\top_\mathcal{I} e_i \}_{i \in I_l(x, \xi)} \} \right).
\end{equation}

Rather than using directions from the approximate normal cone to compute steps, we rely 
on the polar of this cone, called the \emph{approximate tangent cone} and defined by
\begin{equation}
\label{eq:tangentcone}
	T(x,\xi) \; := \; \left\{ v \in \R^n\ 
	| v^\T u \le 0\ \forall u \in N(x,\xi) \right\}.
\end{equation}

An important property of the approximate tangent cone is that it approximates the 
feasible region around $x$, and that moving along all its directions for a 
distance of $\xi$ from $x$ does not break feasibility~\cite{gratton2019direct}. 
Lemma~\ref{lem:feas} below 
provides a formal description of this property (see 
\cite[Lemma~2.1]{gratton2019direct} which is based on \cite[Proposition~2.2]{kolda2007stationarity}).

\begin{lemma}
\label{lem:feas}
Let $x \in \mathcal{F}$ and $\xi > 0$. Then, for any vector $\tilde{d} \in T(x,\xi)$ 
such that $\|\tilde{d}\| \leq \xi$, we have $x + W \tilde{d} \in \mathcal{F}$.
\end{lemma}

Direct-search techniques rely on approximate tangent cones to define 
new feasible points in a way that guarantees convergence to first-order stationarity
\cite{kolda2007stationarity}. 

\section{Full-low evaluation framework with linear constraints}
\label{sec:flec}

In this section, we describe our main algorithmic framework, that belongs to the class 
of \FLE{} algorithms. The main idea behind this technique is the combination of two 
categories of steps. On the one hand, \FullEval{} steps, that are produced at a significant 
cost in terms of function evaluations, are used to yield good performance of the method 
especially in the presence of smoothness. On the other hand, \LowEval{} steps are cheaper 
to compute because they require less evaluations, and are often designed to handle 
the presence of noise and/or non-smoothness in the objective function. We first present 
our general algorithm that combines both types of steps, then dedicate a section to each 
category.

The general mechanism of the \FLE{} approach is described in Algorithm~\ref{alg:flalg}.
In this framework, the iteration type, denoted as \(t_k\), determines whether 
\FullEval{} or \LowEval{} is invoked. The subsequent iteration type~$t_{k+1}$ is decided 
within the invoked function itself, possibly through a user-defined condition. We detail 
the conditions for switching from one iteration type to another in the next sections.

\begin{algorithm}
 {
\caption{\FLE{} framework}
  \label{alg:flalg}

   \textbf{Input}: Initial iterate $x_0 \in \mathcal{F}$, low-eval stepsize $\alpha_0>0$ and iteration type $t_0 = \FullEval{}$.

 \textbf{Output}: The final iterate $x_\infty$.
  \begin{algorithmic}[1]
    \State \textbf{For} $k=0,1,2,\dots$
    \State \hspace{0.25cm} {\bf If} $t_k = \FullEval{}$:
    \State  \hspace{0.5cm} Call $[t_{k+1}, x_{k+1}, \alpha_{k+1} ] = \FullEval{} (x_k, \alpha_k) $ to compute a \FullEval{} step.

    \State \hspace{0.25cm} {\bf Else if} $t_k = \LowEval{}$:
    \State  \hspace{0.5cm} Call $[t_{k+1}, x_{k+1}, \alpha_{k+1} ] = \LowEval{} (x_k, \alpha_k) $ to compute a \LowEval{} step.
  \end{algorithmic}
  }
\end{algorithm}

Apart from requiring feasibility of the initial point, note that 
Algorithm~\ref{alg:flalg} is identical to that of Berahas et al.~\cite{ASBerahas_OSohab_LNVicente_2022}, and that linear constraints are assumed 
to be handled upon computation of a \LowEval{} or a \FullEval{} step. In the next 
sections, we detail our choices for computing those steps.
 
\subsection{Full-eval step based on projections}
\label{subsec:fulleval}

\FullEval{} steps can be implemented by building a model of the objective function around 
the current point and minimizing it to define the next point. A popular approach that 
lies within the derivative-free paradigm consists in computing a finite-difference 
gradient approximation to define a search direction, as well as a stepsize computed 
via line search based on this 
approximation~\cite{ASBerahas_OSohab_LNVicente_2022}. We extend here this approach to the 
linearly constrained setting by considering projections onto the feasible set, a 
popular technique for dealing with linear constraints~\cite{DPBertsekas_2016}.

If the $k$-th iteration of Algorithm~\ref{alg:flalg} is a \FullEval{} step, we define 
a search direction $p_k$ based on an approximate gradient $g_k$ computed through 
finite differences. A natural choice consists in using $p_k=-g_k$. In 
Section~\ref{sec:numsetup}, we detail our practical choices based on quasi-Newton 
formulas. Once $p_k$ has been computed, we then seek candidate steps by considering 
the feasible direction 
\begin{equation}
\label{eq:projpt}
	\bar{x}_k \; = \; \PF{x_k+p_k},
\end{equation}
and performing a line search along the direction $\bar{x}_k-x_k$. 
More precisely, we seek the largest value $\beta \in (0,\bar{\beta}]$ such that 
\begin{equation}
\label{eq:sdccons}
    f\left (x_k + \beta (\bar{x}_k-x_k)\right) 
    \; \le \; 
    f(x_k) + c \beta g_k^\top (\bar{x}_k-x_k).
\end{equation}
where $c \in (0,1)$. We will show in 
Section~\ref{sec:cv} that condition~\eqref{eq:sdccons} is 
satisfied for a sufficiently small value of $\beta$.

Algorithm~\ref{alg:fevconsalg} describes the calculation of a \FullEval{} step for the 
$k$-th iteration of Algorithm~\ref{alg:flalg}. Similarly to the unconstrained case~\cite{ASBerahas_OSohab_LNVicente_2022}, we introduce a switching condition\footnote{In the unconstrained case~\cite{ASBerahas_OSohab_LNVicente_2022}, we have proposed the switching condition $\beta \geq \gamma \rho(\alpha_k)$. Both work for the convergence theory, in the sense of Lemma~\ref{the:nsc}.} that 
controls the norm of the \FullEval{} step. A value $\beta$ is accepted if it satisfies 
the decrease condition~\eqref{eq:sdccons} and
\begin{equation}
\label{eq:key_const}
	\beta \; \geq \; \gamma \alpha_k,
\end{equation}
where $\gamma>0$ is independent of $k$. Condition~\eqref{eq:key_const} guarantees 
that $\beta$ does not go below a certain multiple of~$\alpha_k$, which is the 
stepsize used for computing \LowEval{} steps (see Section~\ref{subsec:loweval}). When 
both~\eqref{eq:sdccons} and~\eqref{eq:key_const} are satisfied, we set $\beta_k=\beta$ 
and define the new iterate as $x_k+ \beta_k \left(\bar{x}_k-x_k\right)$. On the other hand, if 
condition~\eqref{eq:key_const} is violated, the \FullEval{} step is skipped.

\begin{algorithm}[H]
 {
\caption{Constrained \FullEval{} Iteration: Feasible Line Search}
  \label{alg:fevconsalg}

      \textbf{Input}: Iterate $x_k \in \mathcal{F}$ and $\alpha_{k}$. Backtracking global parameters $\bar{\beta} \in (0,1], \gamma > 0$, $\tau \in (0,1)$.

 \textbf{Output}: $t_{k+1}$, $x_{k+1}$, and $\alpha_{k+1}$.

  \begin{algorithmic}[1]
  \State Compute the gradient approximation $g_k$ as well as a search direction 
  $p_k$. Compute $\bar{x}_k$ according to~\eqref{eq:projpt}.
  \State Backtracking line-search: Set $\beta=\bar{\beta}$. 
 
  \State  \textbf{While True} 
  \State \hspace{0.25cm} \textbf{if} \eqref{eq:sdccons} is true or \eqref{eq:key_const} is false, \textbf{break}.
  \State \hspace{0.25cm} Set $\beta = \tau \beta$.
     \State  \textbf{If} (\ref{eq:key_const}) is true,  set $\beta_k = \beta$, 
     $x_{k+1}=x_k+\beta_k(\bar{x}_k-x_k)$, 
     and $t_{k+1} = \FullEval{}$. \textbf{Else},  set
     $x_{k+1}=x_k$ and $t_{k+1} = \LowEval{}$.
     \Statex (The $\LowEval{}$ parameter $\alpha_k$ remains unchanged, $\alpha_{k+1} = \alpha_k$; see Algorithm \ref{alg:pdsf}.)
  \end{algorithmic}
  }
\end{algorithm}

The convexity of the set $\calF$ guarantees that all iterates remain within the feasible set. This is evident when expressing $x_{k+1}$ in the form $(1 - \beta_k) x_k + \beta_k \bar{x}_k$.

\subsection{Low-eval step based on feasible descent cones}
\label{subsec:loweval}

$\LowEval{}$ steps are based on the low evaluation paradigm of probabilistic direct search. This approach can be extended to the linearly constrained case as described in Algorithm \ref{alg:pdsf}.
We suppose that a feasible initial point is provided by the user. At every iteration, 
the algorithm uses a finite number of polling directions to seek a new feasible iterate $x^+$ that reduces the objective function value by a sufficient amount 
\begin{equation}
    \label{eq:ds-sdc}
    f(x^+) \; \leq \; f(x) - \rho(\alpha),
\end{equation}
where $\rho$ is a forcing function classically employed in direct-search methods. The characteristics of $\rho$ are specified in Section \ref{sec:convergence_nonsmooth_const}. 

\begin{algorithm}[H]
 {
\caption{Constrained \LowEval{} Iteration: Feasible Direct Search}
  \label{alg:pdsf}

   \textbf{Input}: Iterate $x_k \in \mathcal{F}$ and stepsize $\alpha_k$.
  Direct-search global parameters $\lambda \geq 1$ and $\theta \in (0,1)$.

   \textbf{Output}: $t_{k+1}$, $x_{k+1}$, and $\alpha_{k+1}$.

  \begin{algorithmic}[1]
  \State Generate a finite set $D_k$ of non-zero polling directions.
  \State \textbf{If} a feasible poll point $x_k + \alpha_k d_k$ is found such that \eqref{eq:ds-sdc} is true for some $d_k\in D_k$, set $x_{k+1} = x_k + \alpha_k d_k$ and $\alpha_{k+1} = \lambda \alpha_k$.
  \State \textbf{Else}, set $x_{k+1} =  x_k$ and $\alpha_{k+1} = \theta \alpha_k$.
  \State Decide \textbf{if} $t_{k+1}=\LowEval$ \textbf{or if} $t_{k+1}=\FullEval$ through a user-defined condition.
  \end{algorithmic}
  }
\end{algorithm}

To ensure the feasibility in Line~2, one can choose directions of the form $W \tilde{d}$, where $\tilde{d} \in T(x_k,\alpha_k)$ with a stepsize less or equal than $\alpha_k$, as shown by Lemma~\ref{lem:feas}.

As Line~4 indicates, the user has the discretion to decide the switching condition from \(\LowEval{}\) to \(\FullEval{}\). The only theoretical requirement is the eventual return to \(\FullEval{}\). An example of such condition could involve restricting the sequence to a predetermined maximum of \(\LowEval{}\) iterations. In our actual implementation, this is achieved by limiting the number of unsuccessful \(\LowEval{}\) attempts to equal the count of backtracking steps executed during the preceding \(\FullEval{}\). This approach ensures a balanced distribution of effort between both types of steps.

\section{Convergence Analysis}
\label{sec:cv}

\subsection{Rate of convergence in the smooth case}
\label{subsec:cvsmooth}

In this section, we analyze the behavior of the class of \FLE{} methods in the smooth case. We show that if the \FullEval{} step generates an infinite sequence of iterates, then the norm of $q(x_k)$ converges to zero with a rate of $1/ \sqrt{k}$.
We now introduce the assumptions needed for the analysis, starting with 
standard boundedness and smoothness requirements.

\begin{assumption}
\label{ass:flow}
    The objective function $f$ is bounded below by 
    $f_\low \in \R$, i.e., $f(x) \ge f_\low$ for all $x \in \R^n$.
\end{assumption}

\begin{assumption}
\label{ass:lip}
    The function $f$ is continuously differentiable and its gradient $\nabla f$ is Lipschitz continuous with 
    constant $L>0$.
\end{assumption}

The next assumptions are related to our approximate gradient and stationary 
measure. For any iterate $x_k$ in $\calF$ computed by Algorithm~\ref{alg:flalg}, 
we define
\begin{equation}
\label{eq:qfunctionsxk}
	q_k = \PF{x_k-\nabla f(x_k)}-x_k, 
	\quad
	q_k^g = \PF{x_k-g_k}-x_k
	\quad \mbox{and} \quad 
	q_k^p = \bar{x}_k - x_k = \PF{x_k+p_k}-x_k,
\end{equation}
In our algorithm, we rely on directions defined using $g_k$. Those should 
be close to the negative of that approximate gradient, in the sense of 
Assumption~\ref{ass:pkgk} below.

\begin{assumption}
\label{ass:pkgk}
	For every iteration $k$, \[
	\frac{(-g_k)^\top q^p_k}{\|q^g_k\|\|q^p_k\|} 
	\; \ge \; \kappa
	\quad \mbox{and} \quad
	u_p \|q_k^p\| \; \le \; \|q_k^g\| \; \le \; U_p \|q_k^p\|.
\]
with $u_p>0$, $U_p > 0$ and $\kappa \in (0, 1]$.
\end{assumption}

When $p_k = -g_k$, Assumption~\ref{ass:pkgk} holds with $\kappa = u_p = U_p = 1$. 
Indeed, the first inequality can be proved using the property of the projection~\cite[Proposition 1.1.4(b)]{DPBertsekas_2016} that implies that
	\[
		(x_k - g_k-\bar{x}_k)^\top (x-\bar{x}_k) \; \le \; 0 
		\qquad \forall x \in \calF.
	\]
 Moreover, using $x=x_k$ in the previous inequality as well as 
 $\bar{x}_k - x_k = q^g_k = q^p_k $ gives
	\[
		g_k^\top (\bar{x}_k-x_k) \le -\|\bar{x}_k-x_k\|^2 
		\quad \Rightarrow \quad 
		- g_k^\top q_k^p \ge  \|q^g_k\|\|q^p_k\|,
	\] 
Finally, in order to relate the control the discrepancy between the true 
criticality measure and its approximation using $g_k$, we require the following 
assumption.

\begin{assumption}
\label{ass:accgk}
The approximate gradient $g_k$ computed at $x_k$ satisfies
\begin{equation} \label{eq:accgk}
    \| \nabla f(x_k) - g_k \| \; \le \; u_g \|q_k^g\|,
\end{equation}
where $u_g \in (0,\kappa(1-c))$ is independent of~$k$.
\end{assumption}
This condition generalizes that in full-low methods for unconstrained 
optimization~\cite[Assumption~3.2]{ASBerahas_OSohab_LNVicente_2022}, albeit 
with a restriction on the constant $u_g$ that becomes superfluous in the 
unconstrained setting. Nevertheless, condition~\eqref{eq:accgk} can 
be guaranteed in a finite number of steps when the gradient is estimated using 
finite differences as shown in Section \ref{subsec:cric}.

We now start our analysis by establishing a lower bound on the stepsize $\beta_k$.

\begin{lemma}
\label{le:boundbetaSF}
	Let Assumptions~\ref{ass:lip}, \ref{ass:pkgk}, and 
    \ref{ass:accgk} hold. If the $k$-th iteration is a successful 
	$\FullEval{}$ iteration, then
	\begin{equation}
	\label{eq:boundbetaSF}
		\beta_k \; \ge \; \beta_{\min} \; := \; \min\left\{ \bar{\beta}, 
		\frac{2\tau(\kappa(1-c) -u_g )u_p}{L}\right\}.
	\end{equation}
\end{lemma}

\begin{proof}
If $\beta_k=\bar{\beta}$ satisfies the decrease condition~\eqref{eq:sdccons}, then~\eqref{eq:boundbetaSF} holds trivially. 
  Therefore, we consider the case where $\beta_k < \bar{\beta}$ and at least one backtracking step was performed. We consider $\beta_k = \tau \beta_k^f$ where $\beta_k^f$ represents the final unsuccessful attempt before satisfying the sufficient decrease condition \eqref{eq:sdccons}. This implies that 
	\begin{equation}
	\label{eq:notsuffdecbndbetaSF}
		c \beta_k^f g_k^\top (\bar{x}_k-x_k)
		\; < \; 
		f(x_k+\beta_k^f(\bar{x}_k-x_k))-f(x_k).
	\end{equation}
 
    Using a Taylor expansion of $f$ around $x_k$ on the right-hand side 
	of~\eqref{eq:notsuffdecbndbetaSF}, we obtain the following inequalities  
	\begin{eqnarray}
            c\beta_k^f g_k^\top (\bar{x}_k-x_k)
            & \leq
		&\beta_k^f\nabla f(x_k)^\top (\bar{x}_k-x_k) 
		+ \frac{L}{2}(\beta_k^f)^2 \|\bar{x}_k-x_k\|^2 \nonumber \\
		c \beta_k^f g_k^\top (\bar{x}_k-x_k)
		&\le 
		&\beta_k^f g_k^\top (\bar{x}_k-x_k) 
		+ \beta_k^f [\nabla f(x_k)-g_k]^\top (\bar{x}_k-x_k) \nonumber \\
		&
		&+ \frac{L}{2}(\beta_k^f)^2\|\bar{x}_k-x_k\|^2 \nonumber \\
		0
		&\le 
		&(1-c)\beta_k^f g_k^\top (\bar{x}_k-x_k) 
		+ \beta_k^f [\nabla f(x_k)-g_k]^\top (\bar{x}_k-x_k) \nonumber \\
		&
		&+ \frac{L}{2}(\beta_k^f)^2\|\bar{x}_k-x_k\|^2.
		\label{eq:taylorbetaSF} 
	\end{eqnarray} 
	Using Assumption \ref{ass:pkgk},
	we have
	\[
		(g_k)^\top q_k^p \; \le \; - \kappa \|q_k^g\| \|q_k^p\| 
		\quad \Leftrightarrow \quad 
		g_k^\top (\bar{x}_k-x_k) \; \le \; - \kappa  \|q_k^g\| \|\bar{x}_k-x_k\|,
	\]
	hence 
	\begin{equation}
	\label{eq:taylorbetaSF1}
		(1-c)\beta_k^f g_k^\top (\bar{x}_k-x_k) 
		\; \le \;
		-(1-c) \kappa \beta_k^f  \|q_k^g\| \|\bar{x}_k-x_k\|.
	\end{equation}
	We now turn to the second term in the right-hand side of~\eqref{eq:taylorbetaSF}. 
	Using Cauchy-Schwarz inequality together with Assumption~\ref{ass:accgk}, 
	we obtain
	\begin{eqnarray*}
		[\nabla f(x_k)-g_k]^\top (\bar{x}_k-x_k) 
		&\le 
		&\|\nabla f(x_k)-g_k\| \|\bar{x}_k-x_k\| \\
		&\le 
		&u_g \|q^g_k\| \|\bar{x}_k-x_k\|.
	\end{eqnarray*}
	Overall, we thus obtain that
	\begin{equation}
	\label{eq:taylorbetaSF2}
		\beta_k^f [\nabla f(x_k)-g_k]^\top (\bar{x}_k-x_k) 
		\; \le \;
		u_g \beta_k^f \|q^g_k\| \|\bar{x}_k-x_k\|.
	\end{equation}
	Putting~\eqref{eq:taylorbetaSF1} and~\eqref{eq:taylorbetaSF2} into 
	\eqref{eq:taylorbetaSF}, we obtain
	\begin{eqnarray*}
		0
		&\le 	
		&-(1-c) \kappa \beta_k^f  \|q_k^g\| \|\bar{x}_k-x_k\|
		+ u_g  \beta_k^f  \|q_k^g\| \|\bar{x}_k-x_k\| \\
		&
		&+ \frac{L}{2}(\beta_k^f)^2 \|\bar{x}_k-x_k\|^2 \\
		0 
		&\le 
		&-(\kappa(1-c) -u_g)\beta_k^f \|q_k^g\| \|\bar{x}_k-x_k\|
		+ \frac{L}{2}(\beta_k^f)^2 \|\bar{x}_k-x_k\|^2.
	\end{eqnarray*}

    Using $\kappa(1-c) -u_g \geq 0$ from Assumption~\ref{ass:accgk} 
    together with Assumption \ref{ass:pkgk}, we show
    \[
    0 \; \leq \; -(\kappa(1-c) -u_g) u_q \beta_k^f \| \bar{x}_k-x_k\|^2
		+ \frac{L}{2}(\beta_k^f)^2 \|\bar{x}_k-x_k\|^2  
    \]
	The latter inequality only holds as long as
	\[
		\beta_k^f \; \ge \; \frac{2(\kappa(1-c) - u_g)u_p}{L}.
	\]
	Since $ \displaystyle \beta_k^f = \beta_k/\tau$, we can conclude that
	\[
		\beta_k \; \ge \; \frac{2\tau(\kappa(1-c) - u_g )u_p}{L}.
	\]
	Combining this result with the case $\beta_k=\bar{\beta}$ gives the desired result.
\end{proof}

We can now establish the main result of the smooth case.

\begin{theorem}
\label{th:cvsmooth}
    Let Assumptions~\ref{ass:flow}--\ref{ass:accgk} hold.
    Let $K \ge 1$ be the first iteration such that $ \|q_{K+1}\| = \|\mathcal{P}_{\calF}\left[x_{K+1}-\nabla f(x_{K+1})\right]
        -x_{K+1}\| \leq \epsilon$. Then, to achieve $\| q_{K+1}\| \leq \epsilon$, Algorithm \ref{alg:flalg} takes at most $n_{SF}^K$ successful \FullEval{} iterations, where
    \begin{equation}
    \label{eq:cvsmooth}
    n_{SF}^K \; \leq \; \ceil{L_1 (f(x_0)-f_{\low}) \epsilon^{-2}},
    \end{equation}
    with $\displaystyle L_1 = \frac{(u_g+1)^2 U_p}{c \kappa \beta_{\min}}$.
       
\end{theorem}

\begin{proof}
    We denote by $\mathcal{I}_{SF}^K$ the set of indices 
    corresponding to successful \FullEval{} iterations.
    Let $k \in \mathcal{I}_{SF}^K$. By definition of such 
    an iteration, the 
    sufficient decrease condition~\eqref{eq:sdccons} is 
    satisfied for $x_{k+1}=\bar{x}_k(\beta_k)$, where 
    $\beta_k$ satisfies~\eqref{eq:boundbetaSF}. Moreover, as 
    shown in the proof of Lemma~\ref{le:boundbetaSF}, we 
    have 
    \[
    	g_k^\top (\bar{x}_k-x_k) \; \le \; - \kappa  \|q_k^g\| \|\bar{x}_k-x_k\|.
    \] 
    Overall, we obtain 
    \begin{eqnarray}
    \label{eq:interres}
        f(x_k) - f(x_{k+1}) 
        &\ge 
        &-c \beta_k g_k^\top (\bar{x}_k-x_k)  \\
        &\ge
        &c \kappa \beta_k  \|q_k^g\| \|\bar{x}_k-x_k\| \nonumber \\
        &\ge
        &\frac{c \kappa \beta_{\min}}{U_p} \|q_k^g\|^2. \nonumber
    \end{eqnarray}
    Meanwhile, using Assumption~\ref{ass:accgk} gives
	\[
		\|q_k\| \; \le \; \|q_k-q_k^g\| + \|q_k^g\| \le (u_g+1)\|q_k^g\|. 
	\]
    Therefore, the decrease achieved at iteration $k$ 
    satisfies
    \begin{equation}
    \label{eq:decreaseSF}
        f(x_k)-f(x_{k+1})
        \; \ge \; 
        \frac{c\kappa \beta_{\min}}{U_p(u_g+1)^2}\|q_k\|^2.
    \end{equation}
    We now consider the changes in function values across all 
    iterations in $\{0,\dots,K-1\}$. Since the iterate does not 
    change on unsuccessful iterations and the function value 
    decreases on successful \LowEval~iterations, we have 
    $f(x_k)-f(x_{k+1}) \ge 0$ for all $k \le K-1$. Combining 
    this observation with Assumption~\ref{ass:flow} 
    and~\eqref{eq:decreaseSF} leads to
    \begin{eqnarray*}
        f(x_0)-f_{\low} 
        &\ge 
        &f(x_0)-f(x_K) \\
        &=
        &\sum_{k=0}^{K-1} f(x_k)-f(x_{k+1}) \\
        &\ge 
        &\sum_{k \in \mathcal{I}_{SF}^K} 
        f(x_k)-f(x_{k+1}) \\
        &\ge 
        &\frac{c \kappa \beta_{\min}}{U_p(u_g+1)^2} 
        \sum_{k \in \mathcal{I}_{SF}^K} \|q_k\|^2 \\
        &>
        &\frac{c \kappa \beta_{\min}}{U_p(u_g+1)^2} 
        n_{SF}^K \epsilon^2.
    \end{eqnarray*}
    Re-arranging the terms and using the assumption that for $\|q_k\| > \epsilon$ for $k \leq K$, lead to the desired conclusion.
\end{proof}

The rate~\eqref{eq:cvsmooth} matches existing result for the unconstrained case~\cite{ASBerahas_OSohab_LNVicente_2022}.
This result primarily addresses the number of successful iterations. However, in the context of DFO, the focus shifts more towards the number of function evaluations. Estimating the upper bound on function evaluations needed to achieve $\|q_{K+1}\| \leq \epsilon$ demands careful consideration of various critical aspects. As outlined in Algorithm \ref{alg:flalg}, an iteration could either be a \FullEval{} iteration, which incurs a cost of up to $n+ \log(\beta_{min}/\bar{\beta}) /\ \log(\tau) +1$ function evaluations, or a \LowEval{} iteration, whose cost is primarily determined by the cardinality of the polling set. While the number of successful \FullEval{} is established, the dynamics between consecutive \LowEval{} iterations and unsuccessful \FullEval{} iterations introduce a layer of complexity that makes it challenging to infer their respective numbers. This difficulty already arises in the unconstrained setting~\cite{ASBerahas_OSohab_LNVicente_2022}, and accurately calculating such a bound falls outside the scope of this convergence analysis.

\subsection{Convergence in the non-smooth case}
\label{sec:convergence_nonsmooth_const}

When the smoothness of the function $f$ is not guaranteed, we rely on the properties of the \LowEval{} steps, and in particular on the sufficient decrease guarantees 
certified by the forcing function. To this end, we make the following 
assumptions.

\begin{assumption}\label{ass:rho}
The function $\rho: \mathbb{R}_{>0} \rightarrow \mathbb{R}_{>0}$ is continuous, positive, non-decreasing, and satisfies $\lim_{\alpha \rightarrow 0^+} \rho(\alpha)/\alpha = 0$.
\end{assumption}

An example of such function is $\rho(\alpha) = \alpha^p$ with $p > 1$.
As in the unconstrained setting~\cite{ASBerahas_OSohab_LNVicente_2022}, we 
require the following assumption on the failure of \FullEval{} iterations.

\begin{assumption}\label{ass:FEfail}
	There exists $\epsilon_g>0$ such that for any $k \in I_{SF}$, where $I_{SF}$ 
	denotes the set of successful \FullEval{} iterations, $\|q_k^g\|>\epsilon_g$.
\end{assumption}

However, we still rely in the analysis on the switching condition~(\ref{eq:key_const}), along with the assumption that the \LowEval{} iterations generate an infinite subsequence of iterates to prove that the direct-search parameter~$\alpha_k$ goes to zero.
This result requires the forcing function to satisfy Assumption~\ref{ass:rho} used in the unconstrained regime.

\begin{lemma}
\label{the:nsc}
Let Assumption~\ref{ass:flow}, \ref{ass:rho}, and \ref{ass:FEfail} hold. Assume that the sequence of iterates~$\{ x_k \}$ is bounded. Then, there exists a point $x_*$ and a subsequence $\mathcal{K} \subset \mathcal{I}_{UL}$ of unsuccessful \LowEval{} iterates for which
\[
\begin{array}{ccc}
     \displaystyle\lim_{k \in \mathcal{K}} x_k \; = \; x_* & \text{ and } & \displaystyle \lim_{k \in \mathcal{K}} \alpha_k \; = \; 0.
\end{array}
\label{eq:xstarconst}
\]
\end{lemma}

\begin{proof}
First, suppose that the set $\mathcal{I}_{SF} \cup \mathcal{I}_{UF} \cup \mathcal{I}_{SL}$ is of infinite cardinality, where $\mathcal{I}_{UF}$ and $\mathcal{I}_{SL}$ are the sets of unsuccessful \FullEval{}~and successful \LowEval{}~iterations, respectively. 
Note that this set represents all iterations~$k$ for which $\alpha_k$ does not 
decrease. 

For all successful \FullEval{} iterations $k \in \mathcal{I}_{SF}$, recall 
that~\eqref{eq:sdccons} holds, i.e.
\[
    f(x_k) - f(x_{k+1})  \ge - c \beta_k g_k^\T (\bar{x}_k - x_k)
    \ge \frac{c \kappa \beta_{k}}{U_p} \|q_k^g\|^2.
\]
Furthermore, the condition~(\ref{eq:key_const}) is satisfied for $\beta = \beta_k$, leading to
\begin{equation} \label{eq:dec1const}
	f(x_k)-f(x_{k+1}) \ge \frac{c \kappa \gamma }{U_p} \alpha_k \|q_k^g\|^2 
	\ge \frac{c \kappa \gamma \epsilon_g^2}{U_p} \alpha_k.
\end{equation}
Meanwhile, successful \LowEval{} iterations $k \in \mathcal{I}_{SL}$ achieve sufficient decrease,
\begin{equation} \label{eq:dec2const}
f(x_k) - f(x_{k+1}) \ge  \rho(\alpha_k).
\end{equation}
Note that in \FullEval{} unsuccessful iterations $k \in \mathcal{I}_{UF}$ neither $x_k$ nor $\alpha_k$ changes.

Hence, given that for unsuccessful \LowEval{} iterations ($\mathcal{I}_{UL}$) the function does not decrease, we can sum from $0$ to $k \in \mathcal{I}_{SF} \cup \mathcal{I}_{UF} \cup \mathcal{I}_{SL}$ the inequalities~(\ref{eq:dec1const}) and~(\ref{eq:dec2const}) to obtain
\begin{eqnarray*}
    f(x_0) - f(x_{k+1}) 
    &\ge 
    &\sum_{k \in \mathcal{I}_{SF}} (f(x_k)-f(x_{k+1})) 
    + \sum_{k \in \mathcal{I}_{SL}} (f(x_k)-f(x_{k+1})) \\
    &\ge 
    & \frac{c \kappa \gamma \epsilon_g^2}{U_p} \sum_{k \in \mathcal{I}_{SF}} \alpha_k + \sum_{k \in \mathcal{I}_{SL}}\rho(\alpha_k).
\end{eqnarray*}
By the boundedness (from below) of~$f$, we conclude that the series are summable, which implies that
$\lim_{k \in \mathcal{I}_{SF}} \alpha_k = 0$ or $\lim_{k \in \mathcal{I}_{SL}} \rho(\alpha_k) = 0$ if one of the sets is infinite. 
Since $\alpha$ remains unchanged during unsuccessful \FullEval{} steps, and under Assumption~\ref{ass:rho}, it follows that $\lim\inf_{k \in \mathcal{I}_{SF} \cup \mathcal{I}_{UF} \cup \mathcal{I}_{SL}} \alpha_k = 0$ (and thus there must be an infinite subsequence of unsuccessful \LowEval{} steps driving $\alpha_k$ to zero).
If both $\mathcal{I}_{SF}$ and $\mathcal{I}_{SL}$ are finite, this implies that $\mathcal{I}_{UF}$ is infinite. In this case, given the mechanism of the algorithm,  there must 
also exist an infinite subsequence of unsuccessful \LowEval{} steps driving $\alpha_k$ to zero.
Similarly, if $\mathcal{I}_{UL}$ is infinite, there must exist an infinite sequence of unsuccessful \LowEval{} steps driving~$\alpha_k$ to zero.

Overall, there must be an infinite sequence of  
unsuccessful \LowEval{} steps driving~$\alpha_k$ to zero.
From the boundedness of the sequence of iterates, one can extract a subsequence~$\mathcal{K}$ of that subsequence satisfying the statement of the lemma.
\end{proof}
\vspace{1ex}

We note that this proof also shows that $\alpha_k$ goes to zero for all $k$.
As in the unconstrained case, convergence results are established using the notion of generalized Clarke-Jahn derivative~\cite{FHClarke_1990} at $x$ along a direction~$d$. 
In Theorem~\ref{th:Clarke-stat}, we show that there exists a limit point which is 
Clarke-Jahn stationary, provided the so-called refining directions are dense in the tangent cone. 

\begin{theorem} \label{th:Clarke-stat}
Let Assumption~\ref{ass:flow}, \ref{ass:rho}, and \ref{ass:FEfail} hold. Assume that the sequence of iterates~$\{ x_k \}$ is bounded.
Let the function $f$ be Lipschitz continuous around the point $x_*$ defined in Lemma~\ref{eq:xstarconst}. Let the set of limit points of
\begin{equation} \label{eq: ref-dirconst}
\left\{ \frac{d_k}{\norm{d_k}}, \; d_k \in D_k, k \in \mathcal{K} \right\}
\end{equation}
be dense in the tangent cone $T(x_*)$, where $\mathcal{K} \subset \mathcal{I}_{UL}$ is given in Lemma~\ref{the:nsc}.

Then, $x_*$ is a Clarke-Jahn stationary point, i.e., $f^{\circ}(x_*;d) \geq 0$ for all normalized $d$ in $T(x^*)$.
\end{theorem}

\begin{proof}
The proof follows standard arguments
in~\cite{CAudet_JEDennis_2002,CAudet_JEDennis_2006,LNVicente_ALCustodio_2012}.
Let $\bar{d}$ be a limit point of~(\ref{eq: ref-dirconst}), identified for a certain subsequence $\mathcal{L} \subseteq \mathcal{K}$. Then, from the basic properties of the generalized Clarke-Jahn derivative, and $k \in \mathcal{L}$,
\begin{align*}
    f^{\circ}(x_*;\bar{d}) &  =  \limsup_{\small \begin{array}{c}
    x_k \to x_*, x_k \in \mathcal{F} \\
     \alpha_k \downarrow 0, x_k + \alpha_k \bar{d} \in \mathcal{F}
\end{array}}  \frac{f(x_k + \alpha_k \bar{d}) - f(x_k)}{\alpha_k}\\
    & \geq \limsup_{\small \begin{array}{c}
    x_k \to x_*, x_k \in \mathcal{F} \\
     \alpha_k \downarrow 0, x_k + \alpha_k d_k \in \mathcal{F}
\end{array}} \left\{ \frac{f(x_k + \alpha_k d_k) - f(x_k)}{\alpha_k} - L_f^* \| d_k - \bar{d} \| \right\} \\
    & = \limsup_{\small \begin{array}{c}
    x_k \to x_*, x_k \in \mathcal{F} \\
     \alpha_k \downarrow 0, x_k + \alpha_k d_k \in \mathcal{F}
\end{array}} \left\{ \frac{f(x_k + \alpha_k d_k) - f(x_k)}{\alpha_k}
     + \frac{\rho(\alpha_k)}{\alpha_k} \right\},
\end{align*}
where $L_f^*$ is the Lipschitz constant of $f$ around $x_*$.
Since $k \in \mathcal{L}$ are unsuccessful \LowEval{} iterations, it follows that $f(x_k + \alpha_k d_k) > f(x_k) - \rho(\alpha_k)$ which implies that
\[
\limsup_{\small \begin{array}{c}
    x_k \to x_*, x_k \in \mathcal{F} \\
     \alpha_k \downarrow 0, x_k + \alpha_k d_k \in \mathcal{F}
\end{array}} \frac{f(x_k + \alpha_k d_k) - f(x_k) + \rho(\alpha_k)}{\alpha_k} \; \geq \; 0.
\]
From this and Assumption~\ref{ass:rho}, we obtain $f^{\circ}(x_*;\bar{d}) \geq 0$.
Given the continuity of $f^{\circ}(x_*;\cdot)$, one has for any $d \in T(x_*)$
such that $\|d\|= 1$,
$f^{\circ}(x_*;d) = \lim_{\bar{d} \to d} f^{\circ}(x_*;\bar{d}) \geq 0$.
\end{proof}

\subsection{More on the smooth case (use of finite difference gradients)}
\label{subsec:cric}
Let us return to the smooth case to clarify the imposition of Assumption~\ref{ass:accgk}. Such an assumption is related to the satisfaction of the so-called criticality step in DFO trust-region methods~\cite{ARConn_KScheinberg_LNVicente_2009a,ARConn_KScheinberg_LNVicente_2009b} 
based on fully linear models.
In the context of Algorithm~\ref{alg:fevconsalg}, those models correspond to an 
approximate gradient $g_k$ built from finite differences.

The $i$-th component of the forward finite-differences (FD) approximation of the gradient at $x_k$ is defined as
\begin{equation} \label{eq:FD}
[\nabla_{h_k} f(x_k)]_i \; = \; \frac{f(x_k + h_k e_i) - f(x_k)}{h_k}, \quad i=1,\ldots,n,
\end{equation}
where $h_k$ is the finite difference parameter and $e_i \in \R^n$ is the $i$-th canonical vector. Computing such a gradient approximation costs $n$ function evaluations per iteration, and it is implicitly assumed that such evaluations can be made.
By using a Taylor expansion, the error in the finite-differences gradient (in the smooth and noiseless setting) can be shown~\cite{ARConn_KScheinberg_LNVicente_2009b} to satisfy
\begin{equation} \label{eq:FDbound}
\|\nabla f(x_k) - \nabla_{h_k} f(x_k) \| \; \leq \; \frac{1}{2} \sqrt{n} \, L \, h_k.
\end{equation}

It becomes then clear that one way to ensure Assumption~\ref{ass:accgk} in practice, when $g_k = \nabla_{h_k} f(x_k)$, is to enforce $h_k \leq u_g'  \| q_k^{h_k} \|$, where $ q_k^{h_k} = \PF{x_k- \nabla_{h_k} f(x_k)}-x_k$ and $ \displaystyle u_g'= 2 u_g / (\sqrt{n} L)$.
Enforcing such a condition is expensive but can be rigorously done through a criticality-step type argument (see Algorithm~\ref{alg:cric}).

\begin{algorithm}[H]
 {
\caption{Criticality step: Performed if $h_k > u_g' \norm{ q_k^{h_k}}$}
  \label{alg:fastsr1}
 {\bf Input:} $h_k$, ${ q_k^{h_k}}^{(0)} = q_k^{h_k}$,  and $\omega \in (0,1)$. Let $j=0$.

 \textbf{Output}: $ q_k^{h_k} = {q_k^{h_k}}^{(j)}$ and $h_k$.

\begin{algorithmic}[1]
 \State \textbf{While} $h_k > u_g'  \norm{ {q_k^{h_k}}^{(j)}}$ \textbf{Do}
 \State \hspace{0.25cm} Set $j = j + 1$ and $h_k = \omega^j u_g' \norm{{q_k^{h_k}}^{(0)}}$.
    \State \hspace{0.25cm} Compute $\nabla_{h_k} f(x_k)$ using~(\ref{eq:FD}) and set $ {q_k^{h_k}}^{(j)} = \PF{x_k- \nabla_{h_k} f(x_k)}-x_k$
  \end{algorithmic}
  \label{alg:cric}
  }
\end{algorithm}

Proposition~\ref{prop:crit} shows that Algorithm~\ref{alg:cric} terminates in a finite number
of steps.

\begin{proposition}\label{prop:crit}
Let Assumption~\ref{ass:lip} hold. If $\|q_k\|>0$, then Algorithm~\ref{alg:cric} terminates in finitely many iterations by computing $h_k$ such that the condition 
$h_k \leq u_g' \| q_k^{h_k}\|$ is satisfied.
\end{proposition}

\begin{proof}
Let us suppose that the algorithm loops infinitely. Then, for all $j \geq 1$,
using Step~2 and the satisfaction of the while--condition in Step~1,
\begin{equation} \label{eq:crit1}
     \| {q_k^{h_k}}^{(j)} \|\; \leq \;  \omega^j \| {q_k^{h_k}}^{(0)} \|.
\end{equation}
On the other hand, for all $j \geq 1$, the FD bound~(\ref{eq:FDbound}), followed by Step~2, gives us
\begin{equation} \label{eq:crit2}
\| \nabla f(x_k) - \nabla_{h_k} f(x_k)^{(j)} \| \; \leq \; \frac{1}{2} \sqrt{n}  L \,  \omega^j u_g' \| {q_k^{h_k}}^{(0)} \|.
\end{equation}
Hence, using (\ref{eq:crit1})--(\ref{eq:crit2}), we have
\begin{eqnarray*}
   \| q_k \|
   \le \| q_k - {q_k^{h_k}}^{(j)} \| + \| {q_k^{h_k}}^{(j)} \| 
   &\le &\| \nabla f(x_k) - \nabla_{h_k} f(x_k)^{(j)} \| + \| {q_k^{h_k}}^{(j)} \| \\ 
   &\le &\| \nabla f(x_k) - \nabla_{h_k} f(x_k)^{(j)} \| 
   + \omega^j \| {q_k^{h_k}}^{(0)} \| \\
   &\le &\left( \frac{\sqrt{n} L u_g'}{2}+ 1\right) \omega^j \| {q_k^{h_k}}^{(0)} \|,
\end{eqnarray*}
where the second inequality on the first line comes from the [non-expansiveness of 
orthogonal projection. By taking limits (and noting that $\omega \in (0,1)$), we conclude that $q_k = 0$, which yields a contradiction.
\end{proof}

\section{Numerical setup} \label{sec:numsetup}

In this section, we will first present our implementation choices for the \FLE{} linearly constrained method. The complete MATLAB implementation is available on GitHub\footnote{ \url{https://github.com/sohaboumaima/FLE}}. The repository includes all the necessary algorithms and testing scripts. The numerical environment of our experiments is also introduced (other methods/solvers tested, test problems chosen, and performance profiles). The tests were run using MATLAB R2019b on an Asus Zenbook with 16GB of RAM and an Intel Core i7-8565U processor running at 1.80GHz.

\subsection{Practical \FullEval{} implementation}

In this section, we present a detailed discussion of the implementation of the \linebreak $\FLE{}$ algorithm in the linearly constrained case. Building upon the principles used in the unconstrained case, we introduce a direction $p_k$ that leverages second-order information for faster convergence. Specifically, we define $p_k = -W H_k W^\top g_k$, where~$H_k$ represents an approximation of the inverse Hessian using the Broyden-Fletcher-Goldfarb-Shanno (BFGS) quasi-Newton update \citep{CGBroyden_1970,RFletcher_1970,DGoldfarb_1970,DFShanno_1970}, as described in Algorithm \ref{alg:BFGS-FD}. Here, $W\in \R^{n \times (n-m)}$ denotes an orthonormal basis for the null space of matrix $A$. Notably, due to the positive definiteness of $H_k$, it follows that $W H_k W^\top$ is also positive definite.

Using $W H_k W^\top g_k$ instead of $H_k g_k$ offers two significant advantages. Firstly, the resulting value of $x_k +p_k$ automatically satisfies the equality constraints, since
\[
A(x_k - W H_k W^\top g_k) \; = \; b - A W (H_kW^\top g_k) \;= \;b.
\] 
Secondly, using this direction allows us to compute $W^\top g_k$ rather than directly calculating $g_k$, thus reducing the computational cost of finite differences from $n$ to $n-m$ function evaluations. Indeed, the forward finite-differences approximation can be reduced to the null space of the linear equality constraints:
\begin{equation} \label{eq:FDw}
[W^\top g_k]_i \; = \; \frac{f(x_k + h_k w_i) - f(x_k)}{h_k}, \quad \text{for} \quad i=1,\ldots,n-m,
\end{equation}
where $h_k$ is the finite difference parameter, and $w_i \in \R^n$ is the $i$-th column vector of $W$. In the numerical experiments, the parameter $h_k$ is set to the square root of Matlab's machine precision.

Our \FullEval{} line-search iteration is described in Algorithm~\ref{alg:BFGS-FD}, which includes BFGS updates for the inverse Hessian approximation $H_k$ using~(\ref{eq:BFGS}). Here, $j_k$ refers to the previous \FullEval{} iteration, and $s_k$ and $y_k$ are given in~(\ref{eq:s-and-y}). Notably, in the non-convex case, the inner product $s_k^\top y_k$ cannot be ensured to be positive. To maintain the positive definiteness of the matrix $H_k$, we skip the BFGS update if $s_k^\top y_k < \epsilon_{c} \|s_k\| \|y_k\|$, with $\epsilon_{c} \in (0,1)$ being independent of $k$. In our implementation, we use $\epsilon_{c}=10^{-10}$.

The line search follows the backtracking scheme described in Algorithm~\ref{alg:fevconsalg}, using standard values $\bar{\beta} = 1$ and $\tau=0.5$. A key feature of our \FLE{} methodology that led to rigorous results (see the proof of Lemma~\ref{the:nsc}) is to stop the line search once condition~(\ref{eq:key_const}) is violated. In our implementation, we use:
\begin{equation} \label{eq:forcing-imp}
\gamma \; = \; 1, \quad
\rho(\alpha_k) \; = \; \min(\gamma_1, \gamma_2 \alpha_k^2), \quad \text{with} \quad \gamma_1\;=\;\gamma_2\;=\;10^{-5}.
\end{equation}

For $k=0$, we perform a backtracking line search using $p_0=-W W^\top g_0$ (and update $t_1$ and~$x_1$) as in Algorithm~\ref{alg:fevconsalg} (with constants as in Algorithm~\ref{alg:BFGS-FD}). The initialization of $H_0$ is done as follows: If $t_1=\FullEval{}$, then we set $H_0 = (y_0^\top s_0)/(y_0^\top y_0) I$, in an attempt to make the size of $H_0$ similar to that of $\nabla^2 f(x_0)^{-1}$~\cite{JNocedal_SJWright_2006}. However, if $t_1=\LowEval{}$, we set $H_0=I$.

\begin{algorithm}
\caption{\FullEval{} Iteration: BFGS with FD Gradients}
  \label{alg:BFGS-FD}
 \textbf{Input}: Iterate $x_k$ with $k \geq 1$. Information $(x_{j_k},g_{j_k},H_{j_k})$ from the previous \FullEval{} iteration~$j_k$  (if $k>0$).
  Backtracking parameters $\bar{\beta}>0$ and $\tau \in (0,1)$. Other parameters $\epsilon_{c},\gamma,\gamma_1>0$.

\textbf{Output}: $t_{k+1}$ and $(x_{k+1},H_k,g_k)$. Return the number $nb_k$ of backtrack attempts.

 \begin{algorithmic}[1]
  \State Compute the FD gradient $W^\top g_k = W^\top \nabla_{h_k} f(x_k)$ using~(\ref{eq:FDw}).
  \State Set
  \begin{equation} \label{eq:s-and-y}
  s_k = x_{k} - x_{j_k} \quad \mbox{and} \quad y_{k} = g_k - g_{j_k}.
  \end{equation}
  \State {\bf If} $s_k^\top y_k \geq \epsilon_{c} \|s_k\| \|y_k\|$,
  set
\begin{equation} \label{eq:BFGS}
H_{k} \; = \; \left (I - \frac{s_k y_k^\top}{y_k^\top s_k}\right) H_{j_k}  \left(I - \frac{y_k s_k^\top }{y_k^\top s_k}\right) + \frac{s_k s_k^\top }{y_k^\top s_k}.
\end{equation}
 \State {\bf Else}, set $H_k = H_{j_k}$.
  \State Compute the direction $- W H_k W^\top g_k$.
  \State Perform a backtracking line-search and update $t_{k+1}$ and $x_{k+1}$ as in Algorithm~\ref{alg:fevconsalg}.
  \end{algorithmic}
\end{algorithm}

\subsection{\LowEval{} implementation}
\label{subsec:setuplow}

We now elaborate on our implementation of Algorithm~\ref{alg:pdsf}, and more precisely 
on the calculation of the polling sets. Our algorithm uses positive generators of the 
approximate tangent cones described in Section~\ref{subsec:tangcone}. By describing an
approximate tangent cone as a conic hull of a finite set of vectors, we can then use those vectors as (feasible) directions. 

The problem of finding such positive generators from a description of the cone through 
linear inequalities has attracted significant research in computational geometry, and 
is sometimes referred to as the representation conversion 
problem~\cite{matheiss1980survey}. Recent advances in linearly constrained 
optimization have featured off-the-shelf softwares to compute those 
generators~\cite{beck2020convergence}. 
We follow here a popular approach in the direct-search 
community~\cite{lewisetal2007}, that splits the problem of computing positive 
generators in two cases. In the first case, we are able to leverage the description of 
the approximate normal cone through positive generators given by~\eqref{eq:normalcone} 
to directly define that of the approximate tangent cone. 
In the second case, we compute positive generators for a subset of the cone, and
positive generators of the tangent cone are then obtained by considering the union 
of all these sets for all possible subsets of columns that yield a full row rank 
matrix~\cite{CJPrice_IDCoope_2003}. One drawback of this strategy is that it 
leads to combinatorial explosion in the subsets of columns that must be 
considered and the number of positive generators that are obtained. For this 
reason, several implementations~\cite{TGKolda_RMLewis_VTorczon_2003,
lewisetal2007} have relied on the double description method from computational 
geometry~\cite{KFukuda_AProdon_1996}. This technique can significantly reduce 
the number of generators that are used to describe the approximate tangent 
cone, in the minority of cases where it is needed on standard test 
problems~\cite{lewisetal2007}. 

Our implementation is that of a probabilistic variant of the aforementioned 
approach proposed by Gratton et al.~\cite{gratton2019direct}, in which the approximate tangent cone is 
decomposed into a subspace part and a pointed cone part (i.e. a cone that 
does not contain a straight line). Given a set of generators for the 
approximate tangent cone, we can then replace the subset related to the 
subspace by a direction drawn uniformly at random within that subspace and 
its negative, while we can randomly sample a fraction of the other 
generators corresponding to the pointed cone part. Such an approach reduces 
the number of polling directions even further, while being endowed with 
almost-sure convergence guarantees~\cite[Proposition 7.1]{gratton2019direct}. 
Our implementation follows that of the \texttt{dspfd} MATLAB 
code~\cite{gratton2019direct}, that uses its own implementation of the 
double description method.

\subsection{Other solvers}

We compared the numerical performance of our implementation of \FLE{} (denoted \texttt{constFLE}) to four other 
approaches:
(i) a line-search BFGS method based on FD gradients (as if there were only \FullEval{} iterations), referred to as 
\texttt{constBFGS};
(ii) probabilistic direct search (as if there were only \LowEval{} iterations), referred to as \texttt{dspfd};
(iii) a mesh adaptive direct search solver, \NOMAD{};
(iv) a direct search solver, referred to as~\ps{}.

Given the detailed description of \texttt{constFLE}, \texttt{constBFGS}, and \texttt{dspfd} in 
previous sections, we only elaborate on \NOMAD{} and \ps{} below.
\NOMAD{}~\cite{audet2021nomad} is a solver that implements Mesh Adaptive Direct Search 
(MADS)~\citep{CAudet_JEDennis_2006} under general nonlinear constraints.
The polling directions belong to positive spanning sets that asymptotically cover the unit sphere densely. 
In the case of inequality constraints, the user is allowed to choose to handle them via extreme-barrier, 
progressive-barrier \cite{audet2009progressive} or filter approaches. In our experiments, we choose the 
default option which is progressive-barrier, but note that an extreme-barrier approach would provide 
similar conclusions.
The \ps{} function is a MATLAB's built-in function that comes as a part of the global optimization 
TOOLBOX~\cite{matlab}. This is a directional direct search method that progresses by accepting a point as 
the new iterate if it satisfies a \textit{simple decrease} condition. For bounds and linear constraints, 
\ps{} modifies poll points to be feasible at every iteration, meaning to satisfy all bounds and linear 
constraints. We adopted the default settings in the choice of polling set which uses the Generalized 
Pattern Search strategy~\cite{torczon1997convergence}.

\subsection{Classes of problems tested}

Evaluating optimization methods crucially involves assessing their performance across diverse scenarios. In pursuit of this, we perform experiments on smooth, noisy, 
and non-smooth problems. For each category, the test set is classified into three distinct classes, namely bound constrained problems, general linearly constrained problems, and problems with at least one linear inequality constraint. Detailed dimensions and inequality counts for each problem are provided in the Appendix for reference. 

For smooth bound constrained problems, we selected~$41$ instances from the \cuter{} library. The dimensions of these instances range from~$2$ to~$20$, and the number of bounds varies between~$1$ and~$40$. The relevant details are summarized in Table \ref{tab:bound}. In the context of smooth general linearly constrained problems, we consider a comprehensive set of~$76$ \cuter{} problems. Each of these problems involves at least one linear constraint, which is not a bound on the variable. The dimensions vary from~$2$ to~$24$, and in cases where linear inequalities are present, their count ranges from~$1$ to~$2000$. A detailed overview of these general constrained problems can be found in Tables~\ref{tab:linequ} to \ref{tab:linineequ}.

To investigate the behavior of the optimization solvers on noisy functions, we conduct experiments using perturbed versions of the aforementioned problems. Following the approach of \cite{JJMore_SMWild_2009}, the perturbed functions are formulated as $f(x) = \phi(x)(1 + \xi (x))$, where $\phi$ represents the original smooth function. In this case, $\xi(x)$ is a realization of a uniform random variable $U(-\epsilon_f, \epsilon_f)$. These noisy functions provide valuable insights into the robustness of optimization algorithms in practical scenarios.

To perform a comparison on nonsmooth problems, we considered two different test sets. The first test set is built from our smooth benchmark, and is meant to illustrate the behavior of our method in presence of mild nonsmoothness.
To create such problems, we considered problems with both bounds and general linear constraints, and moved either the general linear constraints or the bounds into the objective function. As a result of this transformation, we generated a total of 52 bound constrained problems and 107 problems with general linear constraints, out of which 52 included at least one inequality constraint. Comprehensive details about these problems are presented in Tables \ref{tab:nonsmoothbound} to \ref{tab:nonsmoothlininequ}.
In generating general linearly constrained optimization problems, we adopt a method where we penalize only the first portion of the bound constraints in certain cases. This prevents the outcome from being dominated solely by linear equality or inequality constraints. We denote this category as "$1/2 B$" in the tables for ease of reference.
As an illustrative example, let us consider the transformation of problem LSQFIT. The original problem is formulated as follows:
\begin{equation}
\begin{aligned}
\min_{x,y} \quad & \sum_{i=1}^5 (a_ix+y - b_i)^2\\
\textrm{s.t.} \quad & x + y \leq 0.85\\
  &  x \geq 0,
\end{aligned}
\end{equation}
where $a = [0.1, 0.3, 0.5, 0.7, 0.9]$ and $b =[0.25, 0.3, 0.625, 0.701, 1.0]$.
After the transformation, the problem becomes:
\begin{equation}
\begin{aligned}
\min_{x,y} \quad & \sum_{i=1}^5 (a_ix+y - b_i)^2 + \lambda |x + y - 0.85|\\
\textrm{s.t.} \quad & x \geq 0,
\end{aligned}
\end{equation}
where $\lambda$ represents the penalty parameter.
Our second test set of nonsmooth problems consists in 14 linearly constrained minimax problems as introduced in \cite{lukvsan2000test}. In this set, the non-smoothness is introduced by the max operator, resulting in less structured non-smoothness compared to the previous set. Detailed information about these problems is provided in Table \ref{tab:nonsmoothmax}.

\section{Numerical Results} \label{sec:numres}

We present numerical results under the form of performance profiles in order to gauge optimization solvers' effectiveness. As outlined in~\cite{EDDolan_JJMore_2002}, these profiles provide a mean of assessing the performance of a designated set of solvers~$\mathcal{S}$ across a given set of problems~$\mathcal{P}$. They are a visual tool where the highest curve corresponds to the solver with the best overall performance.
Let~$t_{p,s}>0$ be a performance measure of the solver $s \in \mathcal{S}$ on the problem $p \in \mathcal{P}$, which in our case was set to the number of function evaluations.
The curve for a solver~$s$ is defined as the fraction of problems where the performance ratio is at most $\alpha$,
\[
\rho_s(\alpha) \; = \; \frac{1}{|\mathcal{P}|} \text{size} \left\{ p \in \mathcal{P} : r_{p,s} \leq \alpha \right\},
\]
where the performance ratio~$r_{p,s}$ is defined as
\[
r_{p,s} \; = \; \frac{t_{p,s}}{\min \{ t_{p,s} : s \in \mathcal{S}\}}.
\] The convention $r_{p,s} = +\infty$ is used when a solver $s$ fails to satisfy the convergence test for problem~$p$.
The convergence test used is
\begin{equation}
\label{eq:convprof}
f(x_0) - f(x) \; \geq \; (1 - \tau) (f(x_0) - f_L),
\end{equation}
where $\tau > 0$ is a tolerance, $x_0$ is the starting point for the problem, and $f_L$ is computed for each problem $p \in \mathcal{P}$ as the smallest value of~$f$ obtained by any solver within a given number of function evaluations.

In our experiments, we use $100(n+1)$ as a maximum number of function evaluations which is what is need for~$100$ simplex gradients.
Solvers with the highest values of $\rho_s(1)$ are the most efficient, and those with the highest values of $\rho_s(\alpha)$, for large $\alpha$, are the most robust. 

\subsection{Smooth problems}

\subsubsection*{Bound constrained problems}
\begin{figure}[h]
\vskip -25ex
\hskip -30ex
\centering
    \captionsetup{justification=centering}
\begin{minipage}{.47\textwidth}
  \includegraphics[scale = 0.6]{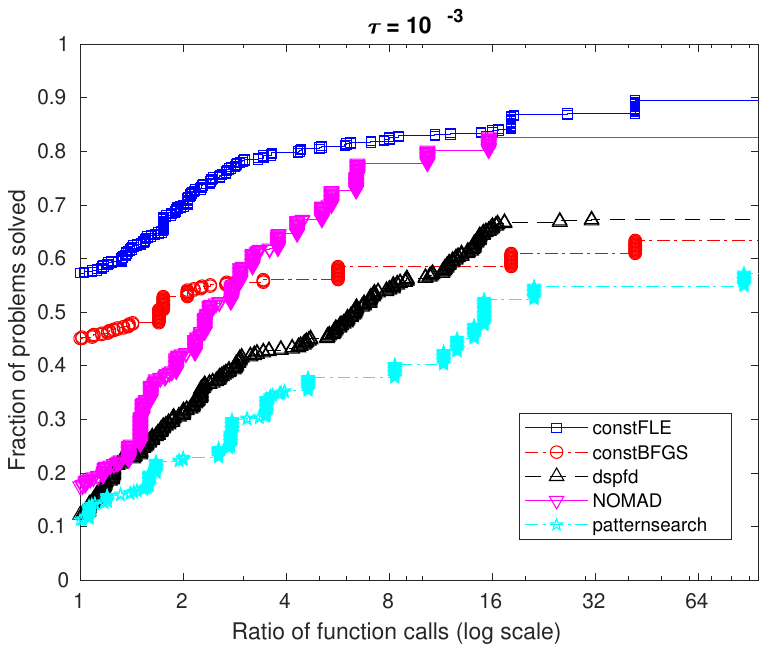}
\end{minipage}%
\hskip 5ex
\begin{minipage}{.47\textwidth}
  \centering
  \includegraphics[scale=0.6]{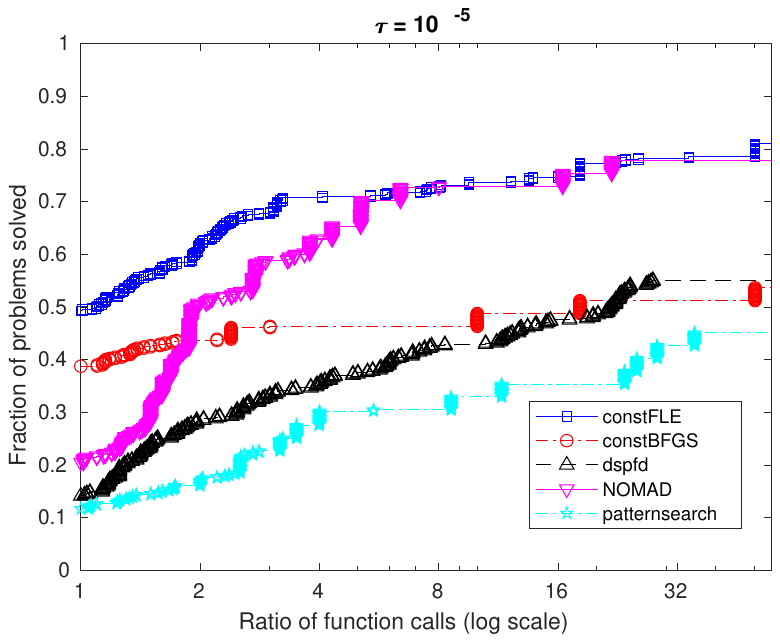}
\end{minipage}
\vskip -25ex
    \caption{Performance profiles with  $\tau = 10^{-3}, 10^{-5}$ of the 5 solvers: constFLE, constBFGS, dspfd,~\NOMAD{}, and~\ps{}. The test set contains 41 smooth bound constrained problems from the~\cuter{} library.}
    \centering
\label{fig:boundsmooth}
\end{figure}

Analyzing those results given in Figure \ref{fig:boundsmooth}, one observes that \FLE{} (blue curve) demonstrates the best performance in terms of efficiency (as indicated by the highest curve at a ratio of 1). It is closely followed by pure \FullEval{} (red curve), with \NOMAD{} (magenta curve) ranking third. When considering robustness, \FLE{} outperforms the others, while \NOMAD{} ranks second. \ps{} ranks last in both efficiency and robustness.
\vspace{-4ex}

\subsubsection*{Linearly constrained problems}
\begin{figure}[h]
\vskip -25ex
\hskip -30ex
\centering
    \captionsetup{justification=centering}
\begin{minipage}{.47\textwidth}
  \includegraphics[scale = 0.6]{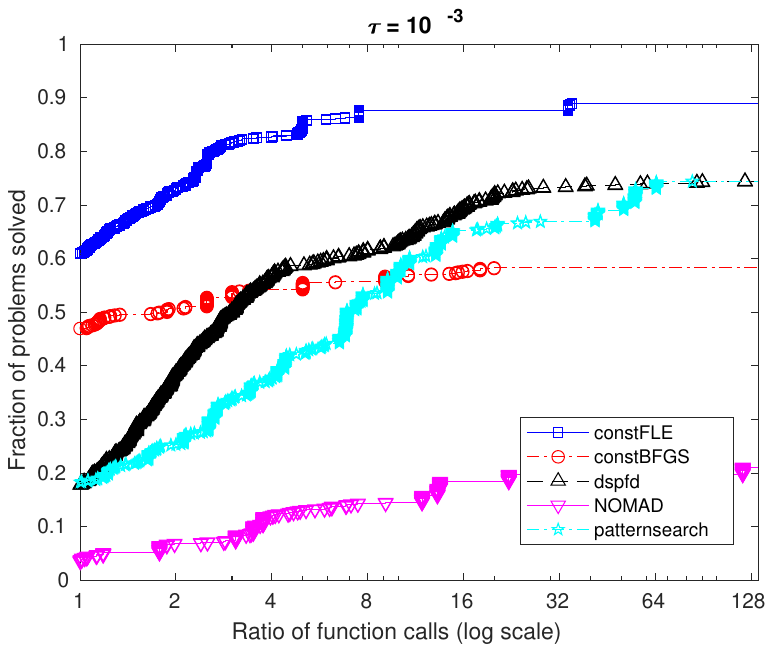}
\end{minipage}%
\hskip 5ex
\begin{minipage}{.47\textwidth}
  \centering
  \includegraphics[scale=0.6]{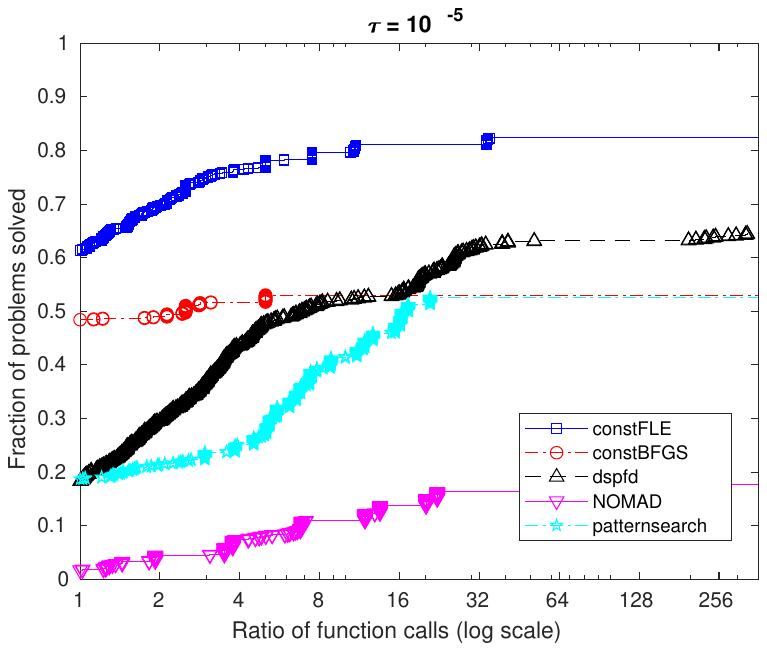}
\end{minipage}
\vskip -25ex
        \caption{Performance profiles with  $\tau = 10^{-3}, 10^{-5}$ of the 5 solvers: constFLE, constBFGS, dspfd,~\NOMAD{}, and~\ps{}. The test set contains 76 smooth problems with general linear constraints from the~\cuter{} library.}
    \centering
\label{fig:gensmooth}
\end{figure}

On general linear equality problems, Figure \ref{fig:gensmooth} illustrates that our method outperforms the four other solvers in terms of both efficiency and robustness. Pure \FullEval{} comes second in terms of efficiency, while pure \LowEval{} performing exceptionally well in terms of robustness and ranks second for that metric. On the other hand,~\ps{} ranks fourth both in terms of efficiency and robustnes, while \NOMAD{} exhibits lower performance due to its limited handling of linear equality constraints, which are present in some of the problems.

\begin{figure}[H]
\vskip -25ex
\hskip -30ex
\centering
    \captionsetup{justification=centering}
\begin{minipage}{.47\textwidth}
  \includegraphics[scale = 0.6]{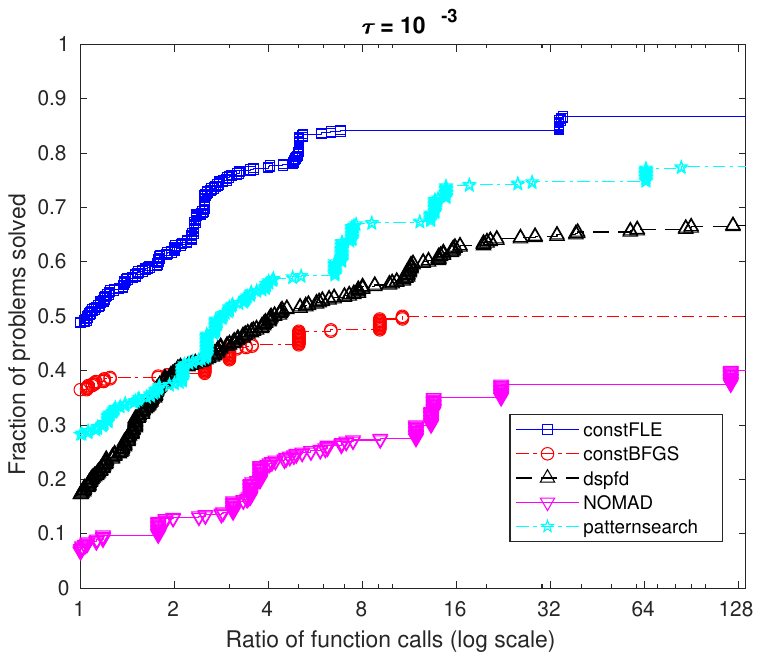}
\end{minipage}%
\hskip 5ex
\begin{minipage}{.47\textwidth}
  \centering
  \includegraphics[scale=0.6]{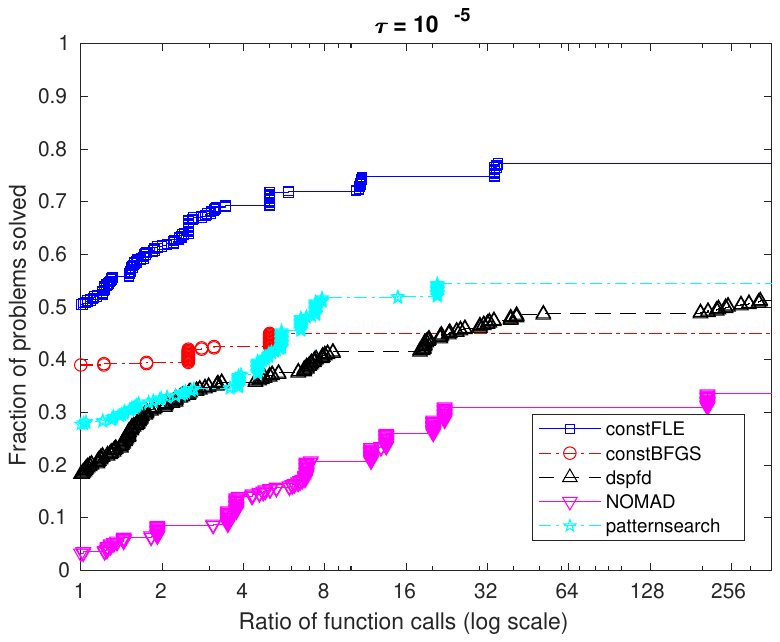}
\end{minipage}
\vskip -25ex
    \caption{Performance profiles with  $\tau = 10^{-3}, 10^{-5}$ of the 5 solvers: constFLE, constBFGS, dspfd,~\NOMAD{}, and~\ps{}. The test set contains 40 smooth problems with at least one inequality constraint from the~\cuter{} library.}
    \centering
\label{fig:linsmooth}
\end{figure}

Figure \ref{fig:linsmooth} provides a more specific comparison of the four solvers on the subset of problems that contain at least one inequality constraint. Even in this context, \FLE{} demonstrates the best performance, followed by \LowEval{}, \FullEval{}, then~\ps{}. It is worth noting that \NOMAD{} shows improved performance compared to the previous experiment given the lack of equality constraints. 

\subsection{Non-smooth problems}

\subsubsection{$\ell_1$ norm problems}

\begin{figure}[h]
\vskip -25ex
\hskip -30ex
\centering
    \captionsetup{justification=centering}
\begin{minipage}{.47\textwidth}
  \includegraphics[scale = 0.6]{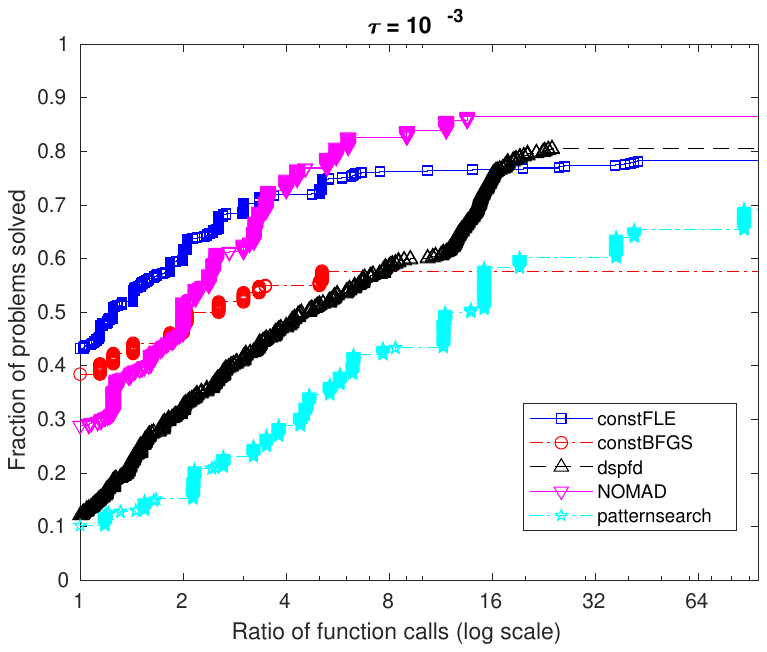}
\end{minipage}%
\hskip 5ex
\begin{minipage}{.47\textwidth}
  \centering
  \includegraphics[scale=0.6]{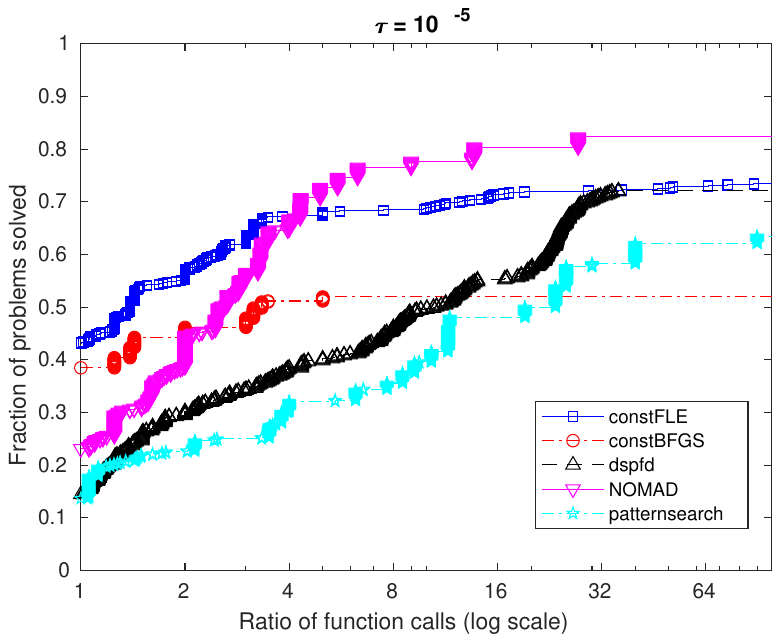}
\end{minipage}
\vskip -25ex
    \caption{Performance profiles with  $\tau = 10^{-3}, 10^{-5}$ of the 5 solvers: constFLE, constBFGS, dspfd,~\NOMAD{}, and~\ps{}. The test set contains 52 non-smooth bound constrained problems.}
    \centering
\label{fig:boundnonsmooth}
\end{figure}

\subparagraph{Bound constrained problems:}
Figure \ref{fig:boundnonsmooth} displays the results obtained from testing non-smooth bound constrained problems. \FLE{} is here the most efficient solver, while \FullEval{} takes the second spot for both low and high accuracy. Meanwhile, \NOMAD{} showcases the best robustness. This observed ranking of solvers in the bound constrained setting remains more or less consistent even with the introduction of the non-smooth regularization.

\newpage

\subparagraph{Linearly constrained problems:}

\begin{figure}[h]
\vskip -25ex
\hskip -30ex
\centering
    \captionsetup{justification=centering}
\begin{minipage}{.47\textwidth}
  \includegraphics[scale = 0.6]{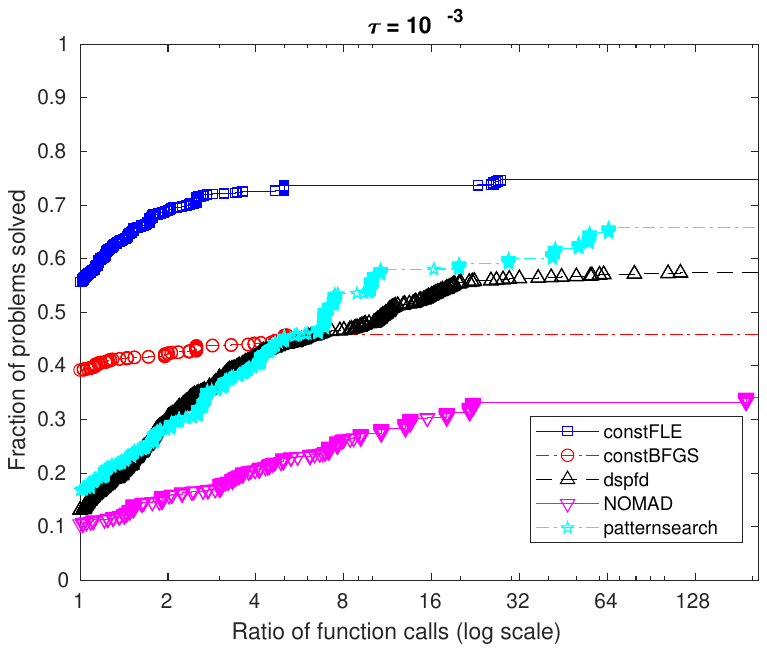}
\end{minipage}%
\hskip 5ex
\begin{minipage}{.47\textwidth}
  \centering
  \includegraphics[scale=0.6]{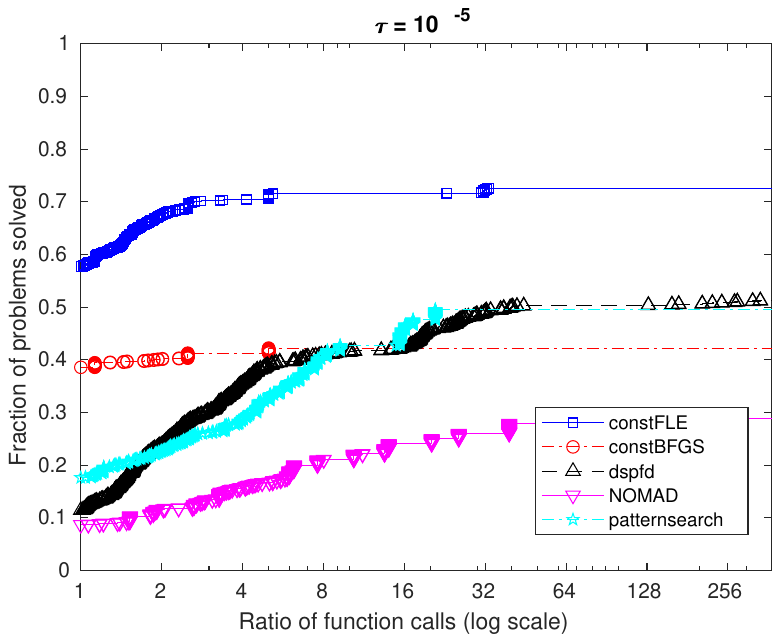}
\end{minipage}
\vskip -25ex
    \caption{Performance profiles with  $\tau = 10^{-3}, 10^{-5}$ of the~5 solvers: constFLE, constBFGS, dspfd,~\NOMAD{}, and~\ps{}. The test set contains 107 non-smooth general linear equality constraints.}
    \centering
\label{fig:gennonsmooth}

\end{figure}

In Figure \ref{fig:gennonsmooth}, we present the results on general non-smooth problems. One can see that the \FLE{} curve is above all, followed by \LowEval{} and \ps{} which exhibit similar performance, \FullEval{}, then \NOMAD{}. As with the smooth case, employing \FLE{} yields better results than using individual steps alone, providing further confirmation of the effectiveness of our approach.

\begin{figure}[h]
\vskip -25ex
\hskip -30ex
\centering
    \captionsetup{justification=centering}
\begin{minipage}{.47\textwidth}
  \includegraphics[scale = 0.6]{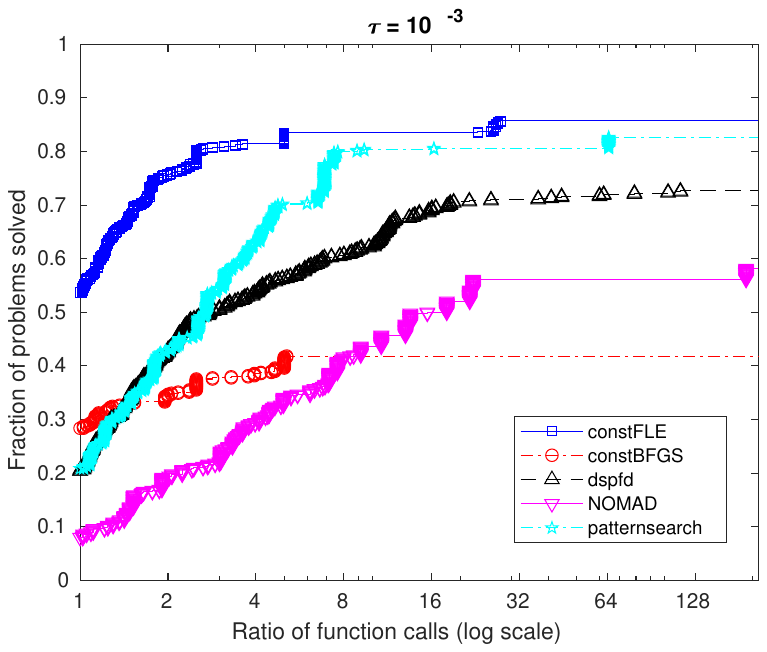}
\end{minipage}%
\hskip 5ex
\begin{minipage}{.47\textwidth}
  \centering
  \includegraphics[scale=0.6]{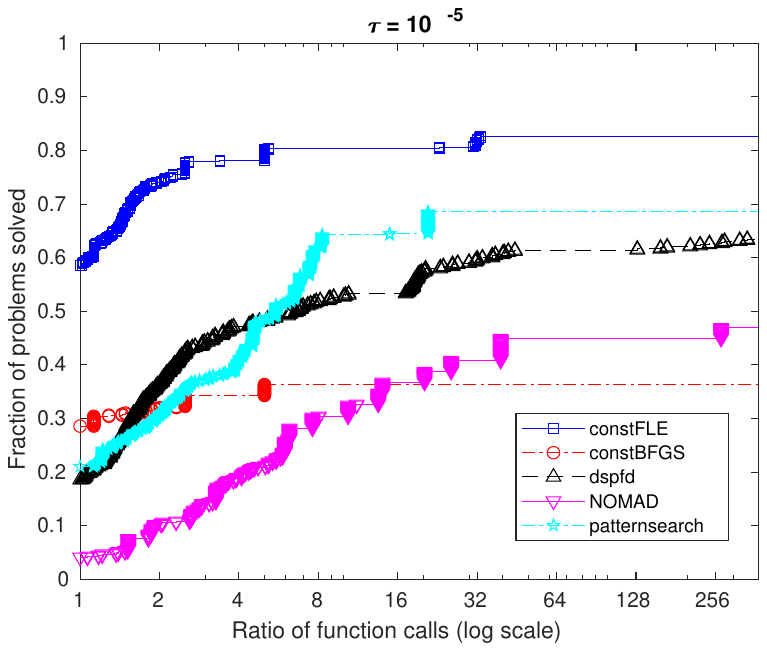}
\end{minipage}
\vskip -25ex
    \caption{Performance profiles with  $\tau = 10^{-3}, 10^{-5}$ of the 5 solvers: constFLE, constBFGS, dspfd,~\NOMAD{}, and~\ps{}. The test set contains 52 non-smooth problems with at least one inequality constraint.}
    \centering
\label{fig:linnonsmooth}
\end{figure}

Furthermore, even within this context, \NOMAD{} faces challenges posed by equality constraints. However, upon their removal as shown in Figure \ref{fig:linnonsmooth}, \NOMAD{} demonstrates improved robustness compared to \FullEval{}.

\subsubsection{Minimax problems}
When the non-smoothness is less structured, methods that estimate the gradient are significantly impacted. Figure \ref{fig:maxnonsmooth} illustrates that the relative performance of these methods is notably different from the earlier observations. On this test set, \NOMAD{} outperforms the other solvers in terms of both efficiency and robustness, followed by \LowEval{}. \FullEval{} comes in third, while \ps{} takes the fourth spot. Among all the methods, \FullEval{} ranks as the least efficient and robust for the reasons aforementioned.

\revisedII{
\begin{figure}[h]
\vskip -25ex
\hskip -30ex
\centering
    \captionsetup{justification=centering}
\begin{minipage}{.47\textwidth}
  \includegraphics[scale = 0.6]{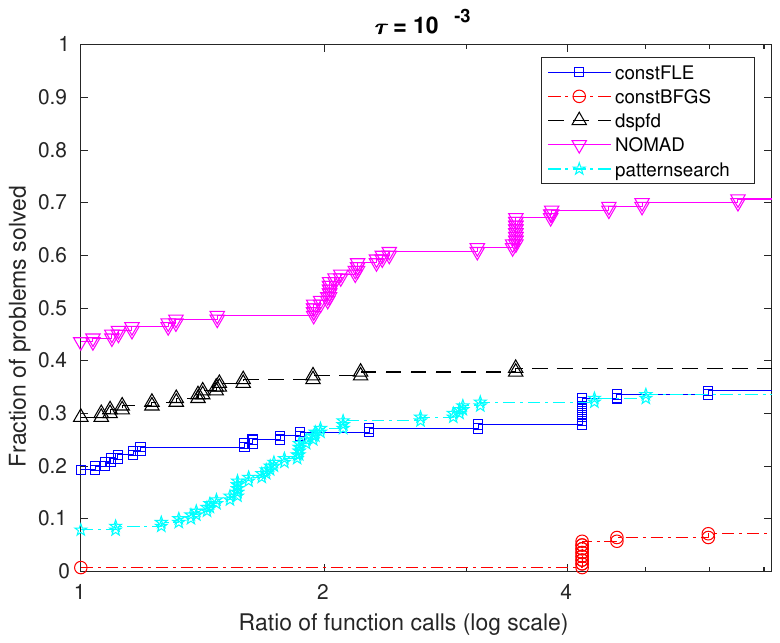}
\end{minipage}%
\hskip 5ex
\begin{minipage}{.47\textwidth}
  \centering
  \includegraphics[scale=0.6]{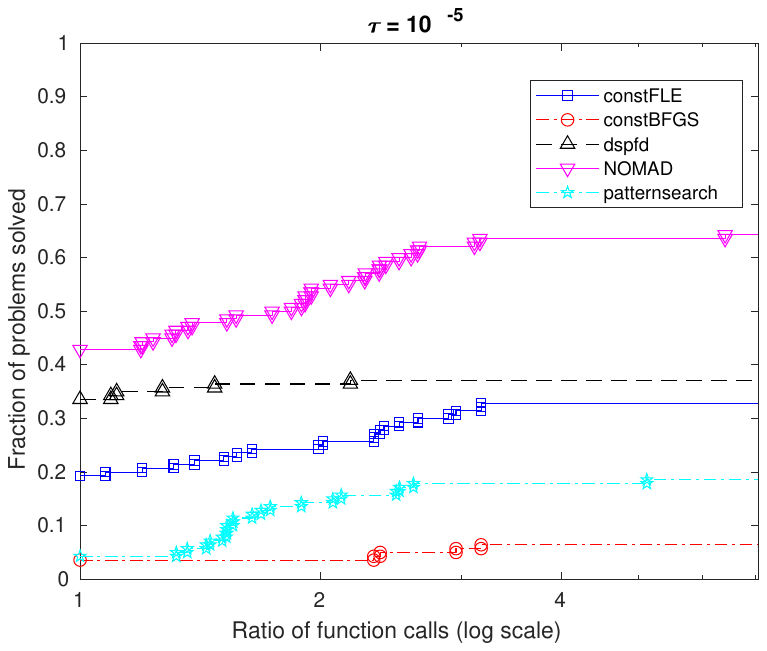}
\end{minipage}
\vskip -25ex
    \caption{Performance profiles with  $\tau = 10^{-3}, 10^{-5}$ of the~5 solvers: constFLE, constBFGS, dspfd,~\NOMAD{}, and~\ps{}. The test set contains 14 bound and linearly constrained problems non-smooth problems.}
    \centering
\label{fig:maxnonsmooth}
\end{figure}}


\subsection{Noisy functions}

\subsubsection*{Bound constrained problems}

\begin{figure}[h]
\vskip -25ex
\hskip -30ex
\centering
    \captionsetup{justification=centering}
\begin{minipage}{.47\textwidth}
  \includegraphics[scale = 0.6]{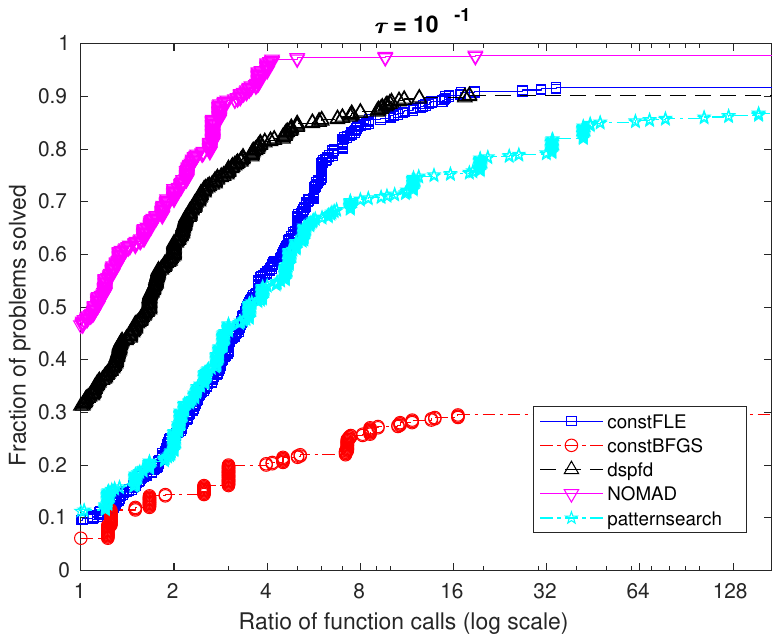}
\end{minipage}%
\hskip 5ex
\begin{minipage}{.47\textwidth}
  \centering
  \includegraphics[scale=0.6]{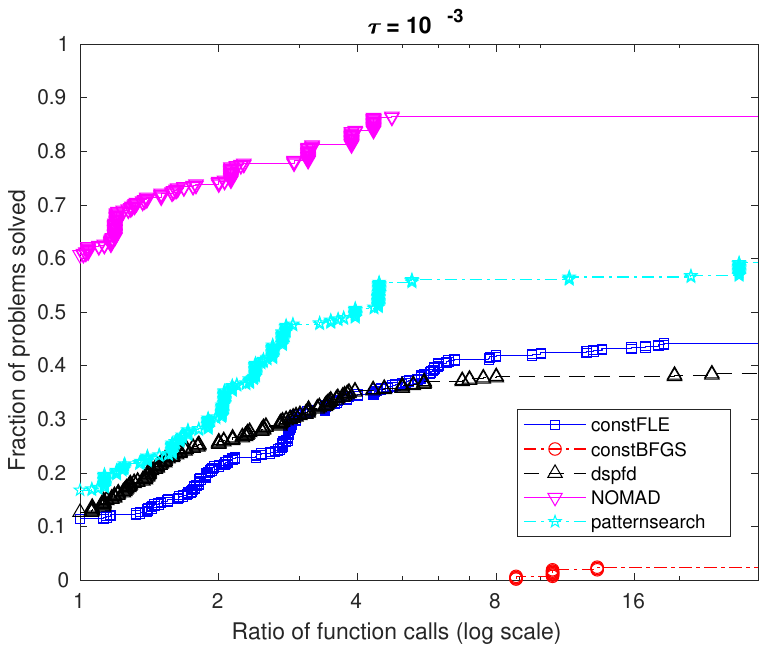}
\end{minipage}
\vskip -25ex
    \caption{Performance profiles with  $\tau = 10^{-1}, 10^{-3}$ of the 5 solvers: constFLE, constBFGS, dspfd,~\NOMAD{}, and~\ps{}. The test set contains 41 noisy bound constrained problems.}
    \centering
\label{fig:boundnoisy}
\end{figure}

In this context, \NOMAD{} demonstrates the best performance in terms of efficiency and robustness. Referring to Figure \ref{fig:boundnoisy}, we can see that the curve corresponding to \FLE{} is between the \FullEval{} and \LowEval{} curves. Such results are conform to observations made in the unconstrained case \cite{ASBerahas_OSohab_LNVicente_2022}, especially for low accuracy. This correspondence arises from \FullEval{} performing poorly when $h$ is equal to the square root of machine precision. Note that \FLE{} is able to outperform \LowEval{} for high accuracy in term of robustness.
On the other hand, at lower accuracies, \ps{} ranks fourth, followed by \FullEval{}. However, its performance improves significantly at higher accuracies, where it ranks second in both efficiency and robustness. This demonstrates its effectiveness in noisy settings.

\subsubsection*{Linearly constrained problems}

\begin{figure}[h]
\vskip -25ex
\hskip -30ex
\centering
    \captionsetup{justification=centering}
\begin{minipage}{.47\textwidth}
  \includegraphics[scale = 0.6]{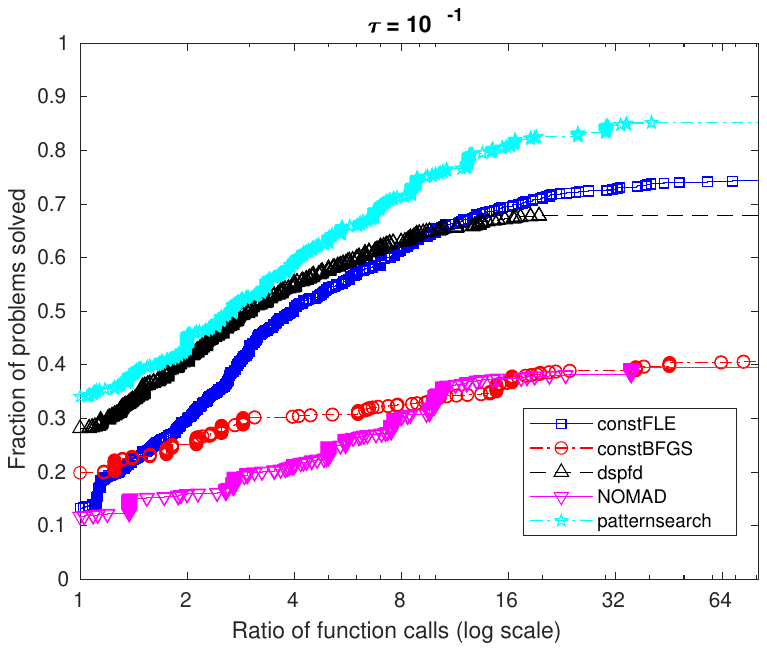}
\end{minipage}%
\hskip 5ex
\begin{minipage}{.47\textwidth}
  \centering
  \includegraphics[scale=0.6]{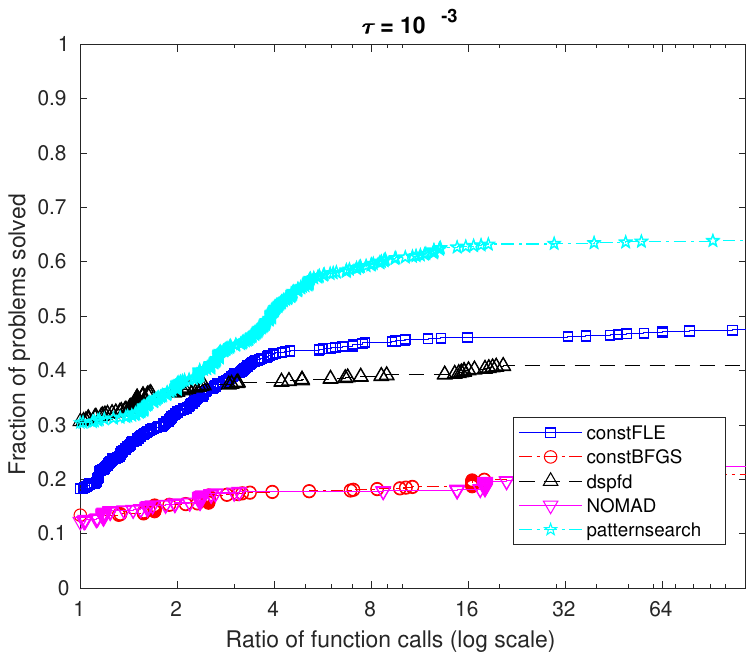}
\end{minipage}
\vskip -25ex
    \caption{Performance profiles with  $\tau = 10^{-3}, 10^{-5}$ of the 5 solvers: constFLE, constBFGS, dspfd,~\NOMAD{}, and~\ps{}. The test set contains 76 noisy problems with general linear constraints.}
    \centering
\label{fig:gennoisy}

\end{figure}

\begin{figure}[h]
\vskip -25ex
\hskip -30ex
\centering
    \captionsetup{justification=centering}
\begin{minipage}{.47\textwidth}
  \includegraphics[scale = 0.6]{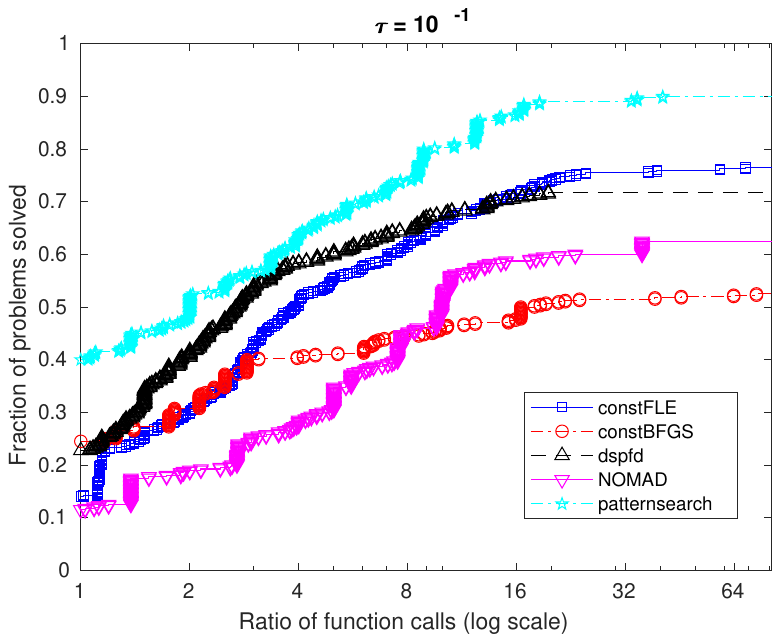}
\end{minipage}%
\hskip 5ex
\begin{minipage}{.47\textwidth}
  \centering
  \includegraphics[scale=0.6]{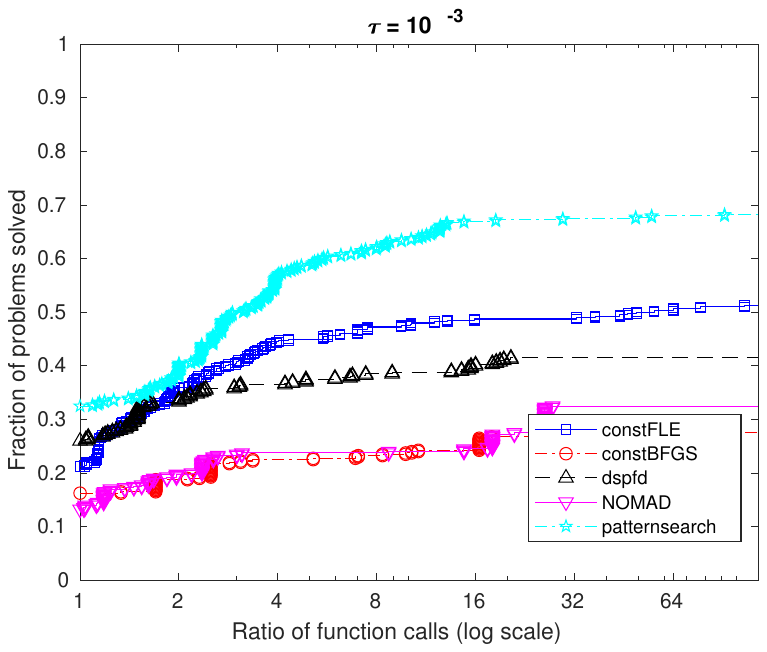}
\end{minipage}
\vskip -25ex
    \caption{Performance profiles with  $\tau = 10^{-3}, 10^{-5}$ of the 5 solvers: constFLE, constBFGS, dspfd,~\NOMAD{}, and~\ps{}. The test set contains 40 noisy problems with at least one inequality constraint.}
    \centering
\label{fig:linnoisy}
\end{figure}

When tested on general linear equality constrained problems,~\ps{} stands out as the most efficient and robust solver. Pure \LowEval{} (probabilistic direct search) shows a comparable efficiency, especially for higher accuracy, followed by \FLE{} which is more robust than~\LowEval{}. Conversely as observed in Figure~\ref{fig:gennoisy}, \NOMAD{} experiences a performance decline similar to observations in both smooth and non-smooth cases. Figure \ref{fig:linnoisy} sheds light on problems featuring linear inequalities. Notably, in this context, \NOMAD{}'s performance stands on par with \FullEval{}, and it even surpasses it, especially under conditions demanding higher accuracy. The relative order of performance among the other solvers remained consistent, with \ps{} demonstrating superior performance.

\section{Conclusions}
\label{sec:conc}

We have proposed an instance of the \FLE{} framework tailored to the presence of bound 
and linear constraints, by combining projected BFGS steps with probabilistic 
direct-search steps within approximate tangent cones. The result method is equipped with 
similar guarantees than in the unconstrained case. In addition, its performance has been 
validated in linearly constrained problems with smooth, non-smooth, and noisy objectives. 
Those experiments overall suggest that our algorithm is able to get the best of both 
worlds, and improve over existing algorithms that do not combine \FullEval{} and 
\LowEval{} steps.

Other variants of the \FLE{} framework may be able to improve on our current 
implementation. In particular, one could rely on trust-region steps as \FullEval{}, 
while one- or two-point feedback feasible approaches that have been proposed more 
generally in the convexly constrained setting. In fact, extending the \FLE{} framework 
to nonlinear, convex constraints is a natural continuation of our work, which may 
benefit from existing results in feedback methods as well as mature theory regarding 
projected gradient techniques.

\section*{Acknowledgments}
This work is partially supported by the U.S. Air Force Office of Scientific Research 
(AFOSR) award FA9550-23-1-0217, and by Agence Nationale de la Recherche through program 
ANR-19-P3IA-0001 (PRAIRIE 3IA Institute).

\section*{Declarations}

{\bf Conflict of interest: }All authors declare that they have no conflict of interest.
\\
{\bf Data Availability: }The data used to support the findings is publicly available. 


\bibliography{fle_const}

\begin{thebibliography}{41}
\providecommand{\natexlab}[1]{#1}
\providecommand{\url}[1]{\texttt{#1}}
\expandafter\ifx\csname urlstyle\endcsname\relax
  \providecommand{\doi}[1]{doi: #1}\else
  \providecommand{\doi}{doi: \begingroup \urlstyle{rm}\Url}\fi

\bibitem[Alarie et~al.(2024)Alarie, Audet, Gheribi, Kokkolaras, and
  Le~Digabel]{alarie2021two}
S.~Alarie, C.~Audet, A.~E. Gheribi, M.~Kokkolaras, and S.~Le~Digabel.
\newblock Two decades of blackbox optimization applications.
\newblock \emph{EURO J. Comput. Optim.}, 9:\penalty0 100011, 2024.

\bibitem[Audet and Dennis~Jr(2009)]{audet2009progressive}
C.~Audet and J.~E. Dennis~Jr.
\newblock A progressive barrier for derivative-free nonlinear programming.
\newblock \emph{SIAM J. Optim.}, 20:\penalty0 445--472, 2009.

\bibitem[Audet and Hare(2017)]{CAudet_WHare_2017}
C.~Audet and W.~Hare.
\newblock \emph{Derivative-Free and Blackbox Optimization}.
\newblock Springer Series in Operations Research and Financial Engineering.
  Springer, Cham, Switzerland, 2017.

\bibitem[Audet and \mbox{Dennis Jr.}(2002)]{CAudet_JEDennis_2002}
C.~Audet and J.~E. \mbox{Dennis Jr.}
\newblock Analysis of generalized pattern searches.
\newblock \emph{SIAM J. Optim.}, 13:\penalty0 889--903, 2002.

\bibitem[Audet and \mbox{Dennis Jr.}(2006)]{CAudet_JEDennis_2006}
C.~Audet and J.~E. \mbox{Dennis Jr.}
\newblock Mesh adaptive direct search algorithms for constrained optimization.
\newblock \emph{SIAM J. Optim.}, 17:\penalty0 188--217, 2006.

\bibitem[Audet et~al.(2021)Audet, Digabel, Montplaisir, and
  Tribes]{audet2021nomad}
C.~Audet, S.~L. Digabel, V.~R. Montplaisir, and C.~Tribes.
\newblock Nomad version 4: Nonlinear optimization with the mads algorithm.
\newblock \emph{arXiv preprint arXiv:2104.11627}, 2021.

\bibitem[Beck and Hallak(2020)]{beck2020convergence}
A.~Beck and N.~Hallak.
\newblock On the convergence to stationary points of deterministic and
  randomized feasible descent directions methods.
\newblock \emph{SIAM J. Optim.}, 30:\penalty0 56--79, 2020.

\bibitem[{Berahas} et~al.(2019){Berahas}, {Byrd}, and
  Nocedal]{ASBerahas_RHByrd_JNocedal_2019}
A.~S. {Berahas}, R.~H. {Byrd}, and J.~Nocedal.
\newblock Derivative-free optimization of noisy functions via quasi-newton
  methods.
\newblock \emph{SIAM J. Optim.}, 29:\penalty0 965--993, 2019.

\bibitem[{Berahas} et~al.(2022){Berahas}, Sohab, and
  {Vicente}]{ASBerahas_OSohab_LNVicente_2022}
A.~S. {Berahas}, O.~Sohab, and L.~N. {Vicente}.
\newblock Full-low evaluation methods for derivative-free optimization.
\newblock \emph{Optim. Methods Softw.}, 38:\penalty0 386--411, 2022.

\bibitem[{Bertsekas}(2016)]{DPBertsekas_2016}
D.~P. {Bertsekas}.
\newblock \emph{Nonlinear Programming}.
\newblock Athena Scientific, Belmont, MA, third edition, 2016.

\bibitem[Broyden(1970)]{CGBroyden_1970}
C.~G. Broyden.
\newblock The convergence of a class of double-rank minimization algorithms 1.
  {G}eneral considerations.
\newblock \emph{IMA J. Appl. Math.}, 6:\penalty0 76--90, 1970.

\bibitem[Clarke(1983)]{FHClarke_1990}
F.~H. Clarke.
\newblock \emph{Optimization and Nonsmooth Analysis}.
\newblock John Wiley \& Sons, New York, 1983.
\newblock Reissued by SIAM, Philadelphia, 1990.

\bibitem[{Conn} et~al.(2000){Conn}, {Gould}, and
  {Toint}]{ARConn_NIMGould_PhLToint_2000}
A.~R. {Conn}, N.~I.~M. {Gould}, and P.~L. {Toint}.
\newblock \emph{Trust-Region Methods}.
\newblock MPS-SIAM Series on Optimization. Society for Industrial and Applied
  Mathematics, Philadelphia, 2000.

\bibitem[{Conn} et~al.(2009){Conn}, Scheinberg, and
  {Vicente}]{ARConn_KScheinberg_LNVicente_2009a}
A.~R. {Conn}, K.~Scheinberg, and L.~N. {Vicente}.
\newblock Global convergence of general derivative-free trust-region algorithms
  to first- and second-order critical points.
\newblock \emph{SIAM J. Optim.}, 20:\penalty0 387--415, 2009.

\bibitem[Conn et~al.(2009)Conn, Scheinberg, and
  Vicente]{ARConn_KScheinberg_LNVicente_2009b}
A.~R. Conn, K.~Scheinberg, and L.~N. Vicente.
\newblock \emph{Introduction to Derivative-Free Optimization}.
\newblock MPS-SIAM Series on Optimization. SIAM, Philadelphia, 2009.

\bibitem[Cust\'odio et~al.(2017)Cust\'odio, Scheinberg, and
  Vicente]{ALCustodio_KScheinberg_LNVicente_2017}
A.~L. Cust\'odio, K.~Scheinberg, and L.~N. Vicente.
\newblock Methodologies and software for derivative-free optimization.
\newblock In T.~Terlaky, M.~F. Anjos, and S.~Ahmed, editors, \emph{Chapter 37
  of Advances and Trends in Optimization with Engineering Applications},
  MOS-SIAM Book Series on Optimization. SIAM, Philadelphia, 2017.

\bibitem[{Dolan} and {Mor\'{e}}(2002)]{EDDolan_JJMore_2002}
E.~D. {Dolan} and J.~J. {Mor\'{e}}.
\newblock Benchmarking optimization software with performance profiles.
\newblock \emph{Math. Program.}, 91:\penalty0 201--213, 2002.

\bibitem[Fletcher(1970)]{RFletcher_1970}
R.~Fletcher.
\newblock A new approach to variable metric algorithms.
\newblock \emph{The Computer Journal}, 13:\penalty0 317--322, 1970.

\bibitem[Fukuda and Prodon(1996)]{KFukuda_AProdon_1996}
K.~Fukuda and A.~Prodon.
\newblock Double description method revisited.
\newblock In M.~Deza, R.~Euler, and I.~Manoussakis, editors,
  \emph{Combinatorics and Computer Science: 8th Franco-Japanese and 4th
  Franco-Chinese Conference, Brest, France, July 3--5, 1995 Selected Papers},
  pages 91--111. Springer, 1996.

\bibitem[Goldfarb(1970)]{DGoldfarb_1970}
D.~Goldfarb.
\newblock A family of variable-metric methods derived by variational means.
\newblock \emph{Math. Comp.}, 24:\penalty0 23--26, 1970.

\bibitem[Gratton et~al.(2011)Gratton, {Toint}, and
  {Tr\"{o}ltzsch}]{SGratton_PhLToint_ATroltzsch_2011}
S.~Gratton, P.~{Toint}, and A.~{Tr\"{o}ltzsch}.
\newblock An active-set trust-region method for derivative-free nonlinear
  bound-constrained optimization.
\newblock \emph{Optim. Methods Softw.}, 21:\penalty0 873--894, 2011.

\bibitem[Gratton et~al.(2019)Gratton, Royer, Vicente, and
  Zhang]{gratton2019direct}
S.~Gratton, C.~W. Royer, L.~N. Vicente, and Z.~Zhang.
\newblock Direct search based on probabilistic feasible descent for bound and
  linearly constrained problems.
\newblock \emph{Comput. Optim. Appl.}, 72:\penalty0 525--559, 2019.

\bibitem[Inc.(Oct 2014)]{matlab}
T.~M. Inc.
\newblock Global optimization toolbox, user's guide.
\newblock \emph{Version 3.3}, Oct 2014.

\bibitem[Jahn(1994)]{jahn1994introduction}
J.~Jahn.
\newblock \emph{Introduction to the Theory of Nonlinear Optimization}.
\newblock Springer Nature, 1994.

\bibitem[Kolda et~al.(2003)Kolda, Lewis, and
  Torczon]{TGKolda_RMLewis_VTorczon_2003}
T.~G. Kolda, R.~M. Lewis, and V.~Torczon.
\newblock Optimization by direct search: {N}ew perspectives on some classical
  and modern methods.
\newblock \emph{SIAM Rev.}, 45:\penalty0 385--482, 2003.

\bibitem[Kolda et~al.(2007)Kolda, Lewis, and Torczon]{kolda2007stationarity}
T.~G. Kolda, R.~M. Lewis, and V.~Torczon.
\newblock Stationarity results for generating set search for linearly
  constrained optimization.
\newblock \emph{SIAM J. Optim.}, 17:\penalty0 943--968, 2007.

\bibitem[Larson et~al.(2019)Larson, Menickelly, and
  Wild]{JLarson_MMenickelly_SWild_2019}
J.~Larson, M.~Menickelly, and S.~Wild.
\newblock Derivative-free optimization methods.
\newblock \emph{Acta Numer.}, 28:\penalty0 287--404, 2019.

\bibitem[{Le Digabel} and {Wild}(2024)]{SLeDigabel_SMWild_2024}
S.~{Le Digabel} and S.~M. {Wild}.
\newblock A taxonomy of constraints in black-box simulation-based optimization.
\newblock \emph{Optim. Eng.}, 25:\penalty0 1125--1143, 2024.

\bibitem[Lewis and Torczon(1999)]{RMLewis_VTorczon_1999}
R.~M. Lewis and V.~Torczon.
\newblock Pattern search algorithms for bound constrained minimization.
\newblock \emph{SIAM J. Optim.}, 9:\penalty0 1082--1099, 1999.

\bibitem[Lewis and Torczon(2000)]{RMLewis_VTorczon_2000}
R.~M. Lewis and V.~Torczon.
\newblock Pattern search methods for linearly constrained minimization.
\newblock \emph{SIAM J. Optim.}, 10:\penalty0 917--941, 2000.

\bibitem[Lewis et~al.(2007)Lewis, Shepherd, and Torczon]{lewisetal2007}
R.~M. Lewis, A.~Shepherd, and V.~Torczon.
\newblock Implementing generating set search methods for linearly constrained
  minimization.
\newblock \emph{SIAM J. Sci. Comput.}, 29:\penalty0 2507--2530, 2007.

\bibitem[Luk{\v{s}}an and Vlcek(2000)]{lukvsan2000test}
L.~Luk{\v{s}}an and J.~Vlcek.
\newblock Test problems for nonsmooth unconstrained and linearly constrained
  optimization.
\newblock Technical Report 798, Institut of Computer Science, Academy of
  Sciences of the Czech Republic, 2000.

\bibitem[Matheiss and Rubin(1980)]{matheiss1980survey}
T.~Matheiss and D.~S. Rubin.
\newblock A survey and comparison of methods for finding all vertices of convex
  polyhedral sets.
\newblock \emph{Math. Oper. Res.}, 5:\penalty0 167--185, 1980.

\bibitem[{Mor\'{e}} and {Wild}(2009)]{JJMore_SMWild_2009}
J.~{Mor\'{e}} and S.~M. {Wild}.
\newblock Benchmarking derivative-free optimization algorithms.
\newblock \emph{SIAM J. Optim.}, 20:\penalty0 172--191, 2009.

\bibitem[Nocedal and Wright(2006)]{JNocedal_SJWright_2006}
J.~Nocedal and S.~J. Wright.
\newblock \emph{Numerical Optimization}.
\newblock Springer-Verlag, Berlin, second edition, 2006.

\bibitem[{Price} and {Coope}(2003)]{CJPrice_IDCoope_2003}
C.~J. {Price} and I.~D. {Coope}.
\newblock Frames and grids in unconstrained and linearly constrained
  optimization: a nonsmooth approach.
\newblock \emph{SIAM J. Optim.}, 14:\penalty0 415--438, 2003.

\bibitem[Rios and Sahinidis(2013)]{LMRios_NVSahinidis_2013}
L.~M. Rios and N.~V. Sahinidis.
\newblock Derivative-free optimization: {A} review of algorithms and comparison
  of software implementations.
\newblock \emph{J. Global Optim.}, 56:\penalty0 1247--1293, 2013.

\bibitem[Shanno(1970)]{DFShanno_1970}
D.~F. Shanno.
\newblock Conditioning of quasi-{N}ewton methods for function minimization.
\newblock \emph{Math. Comp.}, 24:\penalty0 647--656, 1970.

\bibitem[{Shi} et~al.(2023){Shi}, {Xuan}, Oztoprak, and
  Nocedal]{HJMShi_MQXuan_FOztoprak_JNocedal_2023}
H.-J.~M. {Shi}, M.~Q. {Xuan}, F.~Oztoprak, and J.~Nocedal.
\newblock On the numerical performance of finite-difference-based methods for
  derivative-free optimization.
\newblock \emph{Optim. Methods Softw.}, 38:\penalty0 289--311, 2023.

\bibitem[Torczon(1997)]{torczon1997convergence}
V.~Torczon.
\newblock On the convergence of pattern search algorithms.
\newblock \emph{SIAM J. Optim.}, 7:\penalty0 1--25, 1997.

\bibitem[Vicente and Cust\'{o}dio(2012)]{LNVicente_ALCustodio_2012}
L.~N. Vicente and A.~L. Cust\'{o}dio.
\newblock Analysis of direct searches for discontinuous functions.
\newblock \emph{Math. Program.}, 133:\penalty0 299--325, 2012.

\end{thebibliography}

\newpage

\begin{appendices}
\section{List of Problems}
\label{sec:pblist}

The complete list of test problems is presented in Tables \ref{tab:bound} through \ref{tab:nonsmoothmax}. The columns represent various parameters of the problems: Size refers to the dimension of the problem, Bounds indicates the number of bound constraints, LE stands for the number of equality constraints, LI represents the number of inequality constraints, and Func. denotes the number of partial functions in the minimax problem.

\begin{table}[H]

\begin{tabular}{||c c c||} 
 \hline
 Name & Size & Bounds  \\ [0.5ex] 
 \hline\hline
 chenhark & 10 & 10\\ 
 \hline
 explin & 12 & 24  \\
 \hline
 harkerp2 & 10 & 10  \\
 \hline
 hatfldb & 4 & 5  \\
 \hline
 hs3 & 2 & 1  \\ 
 \hline
 hs4 & 2 & 2  \\ 
 \hline
 maxlika & 8 & 16  \\
 \hline
 ncvxbqp1 & 10 & 20  \\ 
 \hline
 oslbqp & 8 & 11  \\ 
 \hline
  pspdoc & 4 & 1  \\ 
 \hline
 weeds & 3 & 4  \\ 
 \hline
 camel6 & 2 & 4  \\ 
 \hline
  eg1 & 3 & 4  \\ 
  \hline
 cvxbqp1 & 10 & 20\\ 
 \hline \hline
 
\end{tabular}
\quad
\begin{tabular}{||c c c||} 
 \hline
 Name & Size & Bounds  \\ [0.5ex] 
 \hline
 \hline
 explin2 & 12 & 24  \\
 \hline
 hart6 & 6 & 12  \\
 \hline
 himmelp1 & 2 & 4  \\
 \hline
 hs2 & 2 & 1  \\ 
 \hline
 hs3mod & 2 & 1  \\ 
 \hline
 hs5 & 2 & 4  \\
 \hline
 mccormck & 10 & 20  \\ 
 \hline
 ncvxbqp2 & 10 & 20  \\ 
 \hline
  palmer1a & 6 & 2  \\ 
 \hline
 palmer4a & 6 & 2  \\
 \hline
 qrtquad & 12 & 12  \\ 
 \hline
 simbqp & 2 & 2  \\ 
 \hline 
 yfit & 3 & 1  \\ 
 \hline
 expquad & 12 & 12\\ 
 \hline \hline
\end{tabular}
\quad
{\renewcommand\arraystretch{1.076}
\begin{tabular}{||c c c||} 
 \hline
 Name & Size & Bounds  \\ [0.5ex] 
 \hline\hline
 hatflda & 4 & 4  \\
 \hline
 hs1 & 2 & 1  \\
 \hline
 hs38 & 4 & 8  \\
 \hline
 hs45 & 5 & 10  \\ 
 \hline
 hs110 & 10 & 20  \\ 
 \hline
 logros & 2 & 2  \\
 \hline
 mdhole & 2 & 1  \\ 
 \hline
 ncvxbqp3 & 10 & 20  \\ 
 \hline
  palmer2b & 4 & 2  \\ 
 \hline
 palmer5b & 9 & 2  \\
 \hline
 probpenl & 10 & 20  \\ 
 \hline
 s368 & 8 & 16  \\ 
 \hline 
 sineali & 20 & 40  \\ 
 \hline \hline
\end{tabular}}
\caption{Bound constrained problems.}
\label{tab:bound}
\end{table}

\begin{table}[H]
\centering
\begin{tabular}{||c c c c||} 
\hline
 Name & Size & Bounds & LE  \\ [0.5ex] 
\hline\hline
aug2d\ & 24 & 0 & 9\\ \hline
bt3 & 5 & 0 & 3\\ \hline
hs28 & 3 & 0 & 1\\ \hline
hs49 & 5 & 0 & 2\\ \hline
hs51 & 5 & 0 & 3\\ \hline
cvxqp2 & 10 & 20 & 2\\ \hline
fccu & 19 & 19 & 8\\ \hline
hs41 & 4 & 8 & 1\\ \hline
hs54 & 6 & 12 & 1\\ \hline
hs62 & 3 & 6 & 1\\ \hline
ncvxqp1 & 10 & 20 & 5\\ \hline
ncvxqp3 & 10 & 20 & 5\\ \hline
ncvxqp5 & 10 & 20 & 2\\ \hline
fits & 10 & 10 & 6\\ \hline
portfl2 & 12 & 24 & 1\\ \hline
portfl4 & 12 & 24 & 1\\ \hline
reading2 & 9 & 14 & 4\\ \hline
sosqp2 & 20 & 40 & 11\\ \hline
\hline
\end{tabular}
\quad
\begin{tabular}{||c c c c||} 
 \hline
 Name & Size & Bounds  & LE\\ [0.5ex] 
 \hline\hline
genhs28 & 10 & 0 & 8\\ \hline
hs9 & 2 & 0 & 1 \\ \hline
hs48 & 5 & 0 & 2\\ \hline
hs50 & 5 & 0 & 3\\ \hline
hs52 & 5 & 0 & 3\\ \hline
cvxqp1 & 10 & 20 & 5\\ \hline
degenlpa & 20 & 40 & 15\\ \hline
hong & 4 & 8 & 1\\ \hline
hs53 & 5 & 10 & 3\\ \hline
hs55 & 6 & 8 & 6\\ \hline
hs112 & 10 & 10 & 3\\ \hline
ncvxqp2 & 10 & 20 & 5\\ \hline
ncvxqp4 & 10 & 20 & 2\\ \hline
ncvxqp6 & 10 & 20 & 2\\ \hline
portfl1 & 12 & 24 & 1\\ \hline
portfl3 & 12 & 24 & 1\\ \hline
portfl6 & 12 & 24 & 1\\ \hline
sosqp1 & 20 & 40 & 11 \\ \hline \hline \end{tabular}
\caption{Linear equality constrained problems.}
\label{tab:linequ}
\end{table}

\begin{table}[H]
\centering
\begin{tabular}{||c c c c c ||} \hline 
Name & Size & Bounds & LE & LI \\ 
[0.5ex] \hline \hline
avgasa &    8 &   16 &    0 &   10 \\ \hline  
biggsc4 &    4 &    8 &    0 &    7 \\ \hline  
dualc2 &    7 &   14 &    1 &  228 \\ \hline  
expfitb &    5 &    0 &    0 &  102 \\ \hline  
hatfldh &    4 &    8 &    0 &    7 \\ \hline  
hs118 &   15 &   30 &    0 &   17 \\ \hline  
hs21mod &    7 &    8 &    0 &    1 \\ \hline  
hs268 &    5 &    0 &    0 &    5 \\ \hline  
hs35mod &    3 &    4 &    0 &    1 \\ \hline  
hs36 &    3 &    6 &    0 &    1 \\ \hline  
hs44 &    4 &    4 &    0 &    6 \\ \hline  
hs76 &    4 &    4 &    0 &    3 \\ \hline  
hs86 &    5 &    5 &    0 &   10 \\ \hline  
lsqfit &    2 &    1 &    0 &    1 \\ \hline  
oet3 &    4 &    0 &    0 & 1002 \\ \hline  
simpllpa &    2 &    2 &    0 &    2 \\ \hline  
sipow1 &    2 &    0 &    0 & 2000 \\ \hline  
sipow2 &    2 &    0 &    0 & 2000 \\ \hline  
sipow3 &    4 &    0 &    0 & 2000 \\ \hline  
stancmin &    3 &    3 &    0 &    2 \\ \hline  
\hline
\end{tabular}
\quad
\begin{tabular}{||c c c c c||} 
 \hline
 Name & Size & Bounds  & LE & LI\\ [0.5ex] 
 \hline\hline
tfi2 &    3 &    0 &    0 &  101 \\ \hline  
avgasb &    8 &   16 &    0 &   10 \\ \hline  
dualc1 &    9 &   18 &    1 &  214 \\ \hline  
dualc5 &    8 &   16 &    1 &  277 \\ \hline  
expfita &    5 &    0 &    0 &   22 \\ \hline  
expfitc &    5 &    0 &    0 &  502 \\ \hline  
hs105 &    8 &   16 &    0 &    1 \\ \hline  
hs21 &    2 &    4 &    0 &    1 \\ \hline  
hs24 &    2 &    2 &    0 &    3 \\ \hline  
hs35 &    3 &    3 &    0 &    1 \\ \hline  
hs37 &    3 &    6 &    0 &    2 \\ \hline  
hs44new &    4 &    4 &    0 &    6 \\ \hline  
hubfit &    2 &    1 &    0 &    1 \\ \hline  
oet1 &    3 &    0 &    0 & 1002 \\ \hline  
pentagon &    6 &    0 &    0 &   15 \\ \hline  
simpllpb &    2 &    2 &    0 &    3 \\ \hline  
sipow1m &    2 &    0 &    0 & 2000 \\ \hline  
sipow2m &    2 &    0 &    0 & 2000 \\ \hline  
sipow4 &    4 &    0 &    0 & 2000 \\ \hline  
zecevic2 &    2 &    4 &    0 &    2 \\ \hline  \hline
\end{tabular}
\caption{Linear inequality constrained problems.}
\label{tab:linineequ}

\end{table}

\begin{table}[H]
\centering
{\renewcommand\arraystretch{1.058}
\begin{tabular}{||p{1.5cm} p{2cm}||}\hline 
Name & Pen. Const \\ [0.5ex] 
\hline \hline
avgasa &  LI \\ \hline  
biggsc4 & LI \\ \hline  
dualc2 &  LE \& LI \\ \hline  
hatfldh &  LI\\ \hline  
hs118 &  LI \\ \hline  
hs21mod &   LI \\ \hline  
hs35mod &   LI \\ \hline  
hs36 &   LI \\ \hline  
hs44 &   LI \\ \hline  
hs76 &   LI \\ \hline  
hs86 &   LI \\ \hline  
lsqfit &    LI \\ \hline  
simpllpa &    LI \\ \hline  
stancmin &    LI \\ \hline  
avgasb &    LI \\ \hline  
dualc1 &    LE \& LI \\ \hline  
dualc5 &    LE \& LI \\ \hline \hline
\end{tabular}}
\quad
{\renewcommand\arraystretch{1.058}
\begin{tabular}{||p{1.5cm} p{2cm}||}
 \hline
 Name & Pen. Const  \\ [0.5ex] 
 \hline\hline
hs105 &    LI \\ \hline  
hs21 &    LI \\ \hline  
hs24 &    LI \\ \hline  
hs35 &    LI \\ \hline  
hs37 &    LI \\ \hline  
hs44new &    LI \\ \hline  
hubfit &    LI \\ \hline  
simpllpb &    LI \\ \hline  
zecevic2 &    LI \\ \hline  
cvxqp2 &   LE \\ \hline  
fccu &  LE \\ \hline  
hs41 &   LE \\ \hline  
hs54 &    LE \\ \hline  
hs62 &    LE \\ \hline  
ncvxqp1 &   LE \\ \hline  
ncvxqp3 &   LE \\ \hline  
ncvxqp5 &   LE \\ \hline  \hline
\end{tabular}}
\quad
{%
\begin{tabular}{||p{1.5cm} p{2cm}||} 
 \hline
 Name & Pen. Const  \\ [0.5ex] 
 \hline\hline
odfits &   LE \\ \hline  
portfl2 &   LE \\ \hline  
portfl4 &   LE \\ \hline  
reading2 &    LE \\ \hline  
sosqp2 &   LE \\ \hline  
cvxqp1 &   LE \\ \hline  
degenlpa &   LE \\ \hline  
hong &    LE \\ \hline  
hs53 &    LE \\ \hline  
hs55 &    LE \\ \hline  
hs112 &   LE \\ \hline  
ncvxqp2 &   LE \\ \hline  
ncvxqp4 &   LE \\ \hline  
ncvxqp6 &   LE \\ \hline  
portfl1 &   LE \\ \hline  
portfl3 &   LE\\ \hline  
portfl6 &   LE \\ \hline  
sosqp1 &   LE \\ \hline  \hline
\end{tabular}}
\caption{Non-smooth bound constrained problems.}
\label{tab:nonsmoothbound}
\end{table}

\begin{table}[H]
\centering
\begin{tabular}{||p{1.5cm} p{2cm}||}\hline 
Name & Pen. Const \\ [0.5ex] 
\hline \hline 
dualc2 &  LI \\ \hline 
dualc1 &  LI \\ \hline  
dualc5 &  LI \\ \hline   
cvxqp2 &   B \\ \hline  
cvxqp2 &   1/2 B \\ \hline  
fccu &  B \\ \hline  
fccu &  1/2 B \\ \hline  
hs41 &   B \\ \hline  
hs41 &   1/2 B \\ \hline  
hs54 &    B \\ \hline
hs54 &    1/2 B\\ \hline
hs62 &    B \\ \hline
hs62 &    1/2 B \\ \hline
ncvxqp1 &   B \\ \hline
ncvxqp1 &   1/2 B \\ \hline  
ncvxqp3 &   B \\ \hline  
ncvxqp3 &   1/2 B \\ \hline 
ncvxqp5 &   B \\ \hline
ncvxqp5 &   1/2 B \\ \hline 
\hline
\end{tabular}
{\renewcommand\arraystretch{1.055}
\quad
\begin{tabular}{||p{1.5cm} p{2cm}||} 
 \hline
 Name & Pen. Const  \\ [0.5ex] 
 \hline\hline 
odfits &   B \\ \hline  
odfits &   1/2 B \\ \hline  
portfl2 &   B \\ \hline
portfl2 &   1/2 B \\ \hline
portfl4 &   B \\ \hline
portfl4 &   1/2 B \\ \hline
reading2 &    B \\ \hline 
reading2 &  1/2 B \\ \hline  
sosqp2 &   B \\  \hline
sosqp2 &   1/2 B \\ \hline
cvxqp1 &   B \\ \hline
cvxqp1 &   1/2 B \\\hline
degenlpa &   B \\\hline
degenlpa &   1/2 B \\\hline
hong &    B \\\hline
hong &    1/2 B \\\hline
hs53 &    B \\\hline
hs53 &    1/2 B \\\hline
\hline
\end{tabular}}
{\renewcommand\arraystretch{1.055}
\quad
\begin{tabular}{||p{1.5cm} p{2cm}||} 
 \hline
 Name & Pen. Const  \\ [0.5ex] 
 \hline\hline 
hs55 &    B \\\hline
hs55 &    1/2 B \\\hline
hs112 &   B \\\hline
hs112 &   1/2 B \\\hline
ncvxqp2 &   B \\\hline
ncvxqp2 &   1/2 B \\\hline
ncvxqp4 &   B \\ \hline
ncvxqp4 &   1/2 B \\ \hline
ncvxqp6 &   B \\ \hline
ncvxqp6 &   1/2 B \\ \hline
portfl1 &   B \\ \hline
portfl1 &   1/2 B \\ \hline
portfl3 &   B \\ \hline
portfl3 &   1/2 B \\ \hline
portfl6 &   B \\ \hline
portfl6 &   1/2 B \\ \hline
sosqp1 &   B \\ \hline
sosqp1 &   1/2 B \\ \hline \hline
\end{tabular}}
\caption{Non-smooth linear equality constrained problems.}
\label{tab:nonsmoothlinequ}
\end{table}

\begin{table}[H]
\centering
{\renewcommand\arraystretch{1.055}
\begin{tabular}{||p{1.5cm} p{2cm}||}
\hline 
Name & Pen. Const \\ [0.5ex] 
\hline \hline 
avgasa &  B \\ \hline 
avgasa &  1/2 B \\ \hline  
biggsc4 &  B \\ \hline   
biggsc4 &   1/2 B \\ \hline  
dualc2 &   B \\ \hline  
dualc2 &  LE \\ \hline  
hatfldh &   B \\ \hline  
hatfldh &    1/2 B \\ \hline
hs118 &    B \\ \hline
hs118 &   1/2 B \\ \hline  
hs21mod &   B \\ \hline
hs21mod &   1/2 B \\ \hline 
hs35mod &    B \\ \hline
hs35mod &    1/2 B \\ \hline 
hs36 &   B \\ \hline 
hs36 &   1/2 B \\ \hline  
hs44 &   B \\ \hline  
\hline
\end{tabular}}
\quad
\begin{tabular}{||p{1.5cm} p{2cm}||} 
 \hline
 Name & Pen. Const  \\ [0.5ex] 
 \hline\hline 
hs44 &    1/2 B \\ \hline  
hs76 &    B \\ \hline  
hs76 &   1/2 B \\ \hline 
hs86 &    B \\ \hline 
hs86 &   1/2 B \\ \hline 
lsqfit &   B \\ \hline 
lsqfit &   1/2 B \\ \hline 
simpllpa &   B \\ \hline  
simpllpa &   LE \\ \hline  
stancmin &   B \\ \hline 
stancmin &    1/2 B \\ \hline 
avgasb &   B \\ \hline  
avgasb &   1/2 B \\ \hline  
dualc1 &   B \\ \hline 
dualc1 &   LE \\ \hline  
dualc5 &    B \\ \hline  
dualc5 &   LE \\ \hline  
hs105 &   B \\ \hline  
\hline
\end{tabular}
{\renewcommand\arraystretch{1.055}
\quad
\begin{tabular}{||p{1.5cm} p{2cm}||} 
 \hline
 Name & Pen. Const  \\ [0.5ex] 
 \hline\hline 
hs105 &   1/2 B \\ \hline 
hs21 &   B \\ \hline  
hs21 &   1/2 B \\ \hline  
hs24 &   B \\ \hline  
hs24 &   1/2 B \\ \hline  
hs35 &    B \\ \hline  
hs35 &    1/2 B \\ \hline  
hs37 &    B \\ \hline  
hs37 &   1/2 B \\ \hline  
hs44new &   B \\ \hline  
hs44new &    1/2 B \\ \hline  
hubfit &   B \\ \hline  
hubfit &   1/2 B \\ \hline  
simpllpb &   B \\ \hline  
simpllpb &   1/2 B \\ \hline  
zecevic2 &   B \\ \hline  
zecevic2 &   1/2 B \\ \hline  \hline
\end{tabular}}
\caption{Non-smooth linear inequality constrained problems.}
\label{tab:nonsmoothlininequ}
\end{table}

\revisedII{\begin{table}[ht]
\centering
\begin{tabular}{||l c c c c c||}
\hline
Name & Size & Func. & Bounds & LE & LI  \\ 
[0.5ex] \hline \hline
MAD1      & 2 & 3 & 0 & 0 & 1 \\ \hline
MAD2      & 2 & 3 & 0 & 0 & 1 \\ \hline
MAD4      & 2 & 3 & 0 & 0 & 1 \\ \hline
MAD5      & 2 & 3 & 0 & 0 & 1 \\ \hline
PENTAGON  & 6 & 3 & 0 & 0 & 15 \\ \hline
MAD6      & 7 & 163 & 1 & 1 & 7 \\ \hline
Wong 2    & 10 & 6 & 0 & 0 & 3 \\ \hline
Wong 3    & 20 & 14 & 0 & 0 & 4 \\ \hline
MAD8     & 20 & 38 & 10 & 0 & 0 \\ \hline
BP filter & 9 & 124 & 0 & 0 & 4 \\ \hline
HS114    & 10 & 9 & 20 & 1 & 4 \\ \hline
Dembo 3  & 7 & 13 & 14 & 0 & 2 \\ \hline
Dembo 5  & 8 & 4 & 16 & 0 & 3 \\ \hline
Dembo 7  & 16 & 19 & 32 & 0 &  1 \\ \hline
\end{tabular}
\caption{Non-smooth minimax linearly constrained problems.}
\label{tab:nonsmoothmax}
\end{table}
}

\end{appendices}
\end{document}